\documentclass[pdflatex,sn-mathphys-num]{sn-jnl}% Math and Physical Sciences Numbered Reference Style
\usepackage{graphicx}%
\usepackage{multirow}%
\usepackage{amsmath,amssymb,amsfonts}%
\usepackage{amsthm}%
\usepackage{mathrsfs}%
\usepackage[title]{appendix}%
\usepackage{xcolor}%
\usepackage{textcomp}%
\usepackage{manyfoot}%
\usepackage{booktabs}%
\usepackage{algpseudocode}%
\usepackage{listings}%
\usepackage{bm}
\usepackage{mathtools}
\usepackage[linesnumbered,boxed,titlenumbered,ruled,nofillcomment]{algorithm2e}

\definecolor{lbcolor}{rgb}{0.9,0.9,0.9}
\theoremstyle{plain}
\newtheorem{theorem}{Theorem} \numberwithin{theorem}{section}

\newtheorem{proposition}[theorem]{Proposition}

\theoremstyle{definition}

\newtheorem{definition}[theorem]{Definition}

\definecolor{mydarkblue}{rgb}{0,0.08,0.45}
\hypersetup{ %
	colorlinks=false,
	linkcolor=mydarkblue,
	citecolor=mydarkblue,
	filecolor=mydarkblue,
	urlcolor=mydarkblue}

\newcommand{\Ibf}[0]{\bm{I}}

\newcommand{\ybf}[0]{\bm{y}}

\newcommand{\Jc}{\mathcal{J}}

\DeclareMathOperator*{\argmin}{arg\,min}

\DeclareMathOperator*{\kr}{\bigodot}

\definecolor{mypink1}{rgb}{0.858, 0.188, 0.478}
\definecolor{mypink2}{RGB}{219, 48, 122}
\definecolor{mypink3}{cmyk}{0, 0.7808, 0.4429, 0.1412}
\definecolor{mygray}{gray}{0.6}

\newlength{\leftstackrelawd}
\newlength{\leftstackrelbwd}
\def\leftstackrel#1#2{\settowidth{\leftstackrelawd}%
{${{}^{#1}}$}\settowidth{\leftstackrelbwd}{$#2$}%
\addtolength{\leftstackrelawd}{-\leftstackrelbwd}%
\leavevmode\ifthenelse{\lengthtest{\leftstackrelawd>0pt}}%
{\kern-.5\leftstackrelawd}{}\mathrel{\mathop{#2}\limits^{#1}}}

 \theoremstyle{plain}

\theoremstyle{definition}

\newtheorem{?}[Th]{Problem}

\newcommand{\bs}[1]{\boldsymbol{#1}}

\renewcommand{\kr}{{\mathrm{kr}}}
\newcommand{\hc}{{\mathrm{HC}}}

\begin{document}

\title[Article Title]{Fast algorithms for 
least square problems with Kronecker lower subsets}


% \author[1]{Osman Asif Malik\thanks{Equal contribution.}}
% \author[2]{Yiming Xu}
% \author[3]{Nuojin Cheng}
% \author[3]{Stephen Becker}
% \author[4]{Alireza Doostan}
% \author[5]{Akil Narayan}

% \affil[1]{Encube Technologies}
% \affil[2]{Department of Mathematics, University of Kentucky }
% \affil[3]{Department of Applied Mathematics, University of Colorado Boulder }
% \affil[4]{Smead Aerospace Engineering Sciences Department, University of Colorado Boulder}
% \affil[5]{Scientific Computing and Imaging Institute, and Department of Mathematics, University of Utah }

\author[1]{\fnm{Osman} \sur{Malik}}\email{}

\author[2]{\fnm{Yiming} \sur{Xu}}\email{}

\author[3]{\fnm{Nuojin} \sur{Cheng}}\email{}

\author[3]{\fnm{Stephen} \sur{Becker}}

\author[4]{\fnm{Alireza} \sur{Doostan}}

\author*[5]{\fnm{Akil} \sur{Narayan}}\email{akil@sci.utah.edu}

\affil[1]{\orgname{Encube Technologies}}

\affil[2]{\orgdiv{Department of Mathematics}, \orgname{University of Kentucky}}

\affil[3]{\orgdiv{Department of Applied Mathematics}, \orgname{University of Colorado Boulder}}

\affil[4]{\orgdiv{Smead Aerospace Engineering Sciences Department}, \orgname{University of Colorado Boulder}}

\affil*[5]{\orgdiv{Scientific Computing and Imaging Institute, and Department of Mathematics}, \orgname{University of Utah}}

% \affil[3]{\orgdiv{Department}, \orgname{Organization}, \orgaddress{\street{Street}, \city{City}, \postcode{610101}, \state{State}, \country{Country}}}


\abstract{
    While leverage score sampling provides powerful tools for approximating solutions to large least squares problems, the cost of computing exact scores and sampling often prohibits practical application. This paper addresses this challenge by developing a new and efficient algorithm for exact leverage score sampling applicable to matrices that are lower column subsets of Kronecker product matrices. We synthesize relevant approximation guarantees and detail the algorithm that specifically leverages this structural property for computational efficiency. Through numerical examples, we demonstrate that utilizing efficiently computed exact leverage scores via our methods significantly reduces approximation errors, as compared to established approximate leverage score sampling strategies when applied to this important class of structured matrices.
}

\keywords{Leverage score sampling, Kronecker product, Sketching, Least squares, Uncertainty quantification}



\maketitle

\section{Introduction}
\label{sec:intro}

Developing efficient algorithms to compute approximate solutions for large-scale linear least squares problems is a well-studied topic. Given a design matrix $\bm{A} \in \mathbb{R}^{M \times N}$ and an observation vector $\bm{b} \in \mathbb{R}^M$, the full least squares problem seeks the optimal solution $\bm{x}^\ast$ that minimizes the residual norm:
\begin{equation}\label{eq:exact-solution}
  \bm{x}^\ast = \argmin_{\bm{x} \in \mathbb{R}^N} \left\| \bm{A} \bm{x} - \bm{b} \right\|_2.
\end{equation}
Approximate solution algorithms are crucial when dealing with large datasets, either because $M$ is excessively large (making direct manipulation or storage of $\bm{A}$ infeasible) or because acquiring all entries of $\bm{b}$ is prohibitively expensive. A standard and widely used approach involves constructing a sketching operator $\bm{S} \in \mathbb{R}^{K \times M}$, typically with $K \ll M$, and solving the smaller, sketched least squares problem:
\begin{equation}\label{eq:sampled-solution}
  \tilde{\bm{x}} = \argmin_{\bm{x} \in \mathbb{R}^N} \left\| \bm{S} \bm{A} \bm{x} - \bm{S} \bm{b} \right\|_2.
\end{equation}
When $\bm{S}$ is generated randomly, the approximate solution $\tilde{\bm{x}}$ often comes with probabilistic error bounds relative to $\bm{x}^\ast$. One prominent technique chooses $\bm{S}$ as a random row sampling or reweighting operator based on the statistical leverage scores of $\bm{A}$ \citep{drineas2006SamplingAlgorithms, drineas2008RelativeerrorCUR, drineas2011FasterLeast, larsen2020PracticalLeverageBased}. This approach typically yields guarantees of the following form: for a target accuracy $\epsilon$ and failure probability $\delta$, choosing $K \gtrsim \frac{N \log N}{\epsilon^2}$ rows appropriately ensures
\begin{equation} \label{eq:leverage_guarantee}
  % K \gtrsim \frac{N \log N}{\epsilon^2} \quad \Longrightarrow \quad
  \| \bm{A} \tilde{\bm{x}} - \bm{b} \|_2^2 \leq (1 + \epsilon) \| \bm{A} \bm{x}^\ast - \bm{b} \|_2^2,
\end{equation}
with probability at least $1-\delta$. A major practical challenge, however, lies in efficiently computing or sampling according to the leverage scores of $\bm{A}$. Direct computation generally requires forming an orthonormal basis for the column space (range) of $\bm{A}$, often via the QR decomposition or SVD, which has a computational cost scaling 
at least linearly with $M$
%poorly with $M$ 
(e.g., $\mathcal{O}(MN^2)$ for QR when $N\le M$). Two main strategies are typically employed to mitigate this cost:
\begin{enumerate}
    \item[(i)] Sample from efficiently computable \textit{approximate} leverage scores, sacrificing exactness for speed \citep{drineas2012FastApproximation}; or, 
    \item[(ii)] Compute or sample from \textit{exact} leverage scores efficiently by exploiting specific structures within $\bm{A}$ using specially developed algorithms.
\end{enumerate}

This paper focuses on the second strategy. We develop efficient algorithms for exact leverage score sampling when the matrix $\bm{A}$ possesses a specific structure, namely, when it is a lower column subset of a Kronecker product matrix. Specifically, consider $D$ factor matrices $\bm{A}^{(1)}, \ldots, \bm{A}^{(D)}$, where $\bm{A}^{(d)} \in \mathbb{R}^{M_d \times N_d}$. Let the full Kronecker product be $\bm{A}_{\textup{kr}} \coloneqq \bigotimes_{d=1}^D \bm{A}^{(d)}$ which has $\prod_{d=1}^D M_d$ rows. The matrices $\bm{A}$ we consider are formed by selecting columns of $\bm{A}_{\textup{kr}}$ whose multi-indices belong to a set $\mathcal{J}$ that satisfies the \textit{lower subset} property, as formally defined in Definition~\ref{def:lower-set}.
This particular matrix structure is motivated by applications in high-dimensional function approximation, particularly when constructing emulators or surrogate models from tensor product function spaces, such as multivariate polynomials \citep{adcock_sparse_2022}, which we will explain in detail later.

\subsection{Contributions}
\label{ssec:contributions}

The primary contributions of this paper are twofold:
\begin{itemize}
    \item We present an algorithm for \emph{exact} leverage score sampling specifically designed for matrices $\bm{A}$ that are lower column subsets of Kronecker product matrices. The algorithm achieves sampling efficiency significantly greater than the standard implementation. To the best of our knowledge, this is the first efficient algorithm proposed for exact leverage score sampling for this specific class of structured matrices.
    \item We provide an implementation of our proposed algorithm and compare its performance experimentally against relevant competing methods across various benchmark problems. The numerical results demonstrate the computational efficiency and competitive accuracy of our approach.
\end{itemize}

\subsection{Related Work}
\label{ssec:related_work}

Several previous works have developed randomized sketching operators $\bm{S}$ that can be applied particularly efficiently when the matrix $\bm{A}$ or its constituent vectors possess Kronecker product structure. Often, efficiency gains are realized when applying the sketch to vectors $\bm{v}$ that are themselves Kronecker products, i.e., vectors of the form
\begin{equation} \label{eq:kronecker-vector}
	\bm{v} = \bm{a}_1 \otimes \bm{a}_2 \otimes \cdots \otimes \bm{a}_D,
\end{equation}
where each $\bm{a}_d \in \mathbb{R}^{M_d}$. We briefly review major classes of such methods below. For a more comprehensive overview, we refer the reader to~\citep[Section~7]{murray2023RandomizedNumerical}.

\paragraph{Row-structured sketches}
A first class of sketches imposes specific structure on the rows of $\bm{S}$, enabling efficient application to vectors of the form \eqref{eq:kronecker-vector}. Biagioni et al.~\citep{biagioni2015RandomizedInterpolative} proposed a sketch where each row is the Kronecker product of $D$ dense random vectors. Sun et al.~\citep{sun2018TensorRandom} independently proposed sketches with identical structure, including a variance-reduced version. These concepts were extended by constructing sketches whose rows correspond to vectorized low-rank tensors, such as those in Canonical Polyadic (CP) or Tensor Train (TT) formats, where the factors or cores contain appropriate random entries \citep{rakhshan2020TensorizedRandomProjectionsAISTATS, rakhshan2021RademacherRandom}. Iwen et al.~\citep{iwen2021LowerMemoryOblivious} introduced a two-stage sketch involving an initial Kronecker product sketch followed by a secondary sketch for further dimension reduction. Additional details on these methods can be found in the review by Martinsson and Tropp~\citep[Sec.~9.4]{martinsson2020RandomizedNumerical}.

\paragraph{Kronecker fast Johnson--Lindenstrauss transforms (KFJLTs)}
A second class enhances the standard fast Johnson--Lindenstrauss transform (FJLT) \citep{ailon2009FastJohnson} with additional structure compatible with Kronecker products, leading to further computational speedups when applied to vectors like \eqref{eq:kronecker-vector}. These are often termed Kronecker FJLTs (KFJLTs) or, when based on the Hadamard transform, Tensor Subsampled Randomized Hadamard Transforms (TensorSRHTs). KFJLTs were introduced by Battaglino et al.~\citep{battaglino2018PracticalRandomized} for tensor decomposition algorithms, with subsequent theoretical and empirical studies in \citep{jin2020FasterJohnsonLindenstrauss, malik2020GuaranteesKronecker, bamberger2021JohnsonLindenstraussEmbeddings}.

\paragraph{TensorSketch}
The TensorSketch represents a third class, viewable as a structured variant of the CountSketch \citep{clarkson2017LowRankApproximation}, specifically designed for efficient application to Kronecker product vectors \eqref{eq:kronecker-vector}. Its development occurred across several works \citep{pagh2013CompressedMatrix, pham2013FastScalable, avron2014SubspaceEmbeddings, diao2018SketchingKronecker}.

\paragraph{Recursive sketches}
Fourth, recursive sketching techniques apply sketches hierarchically to achieve improved theoretical guarantees or efficiency. Ahle et al.~\citep{ahle2020ObliviousSketchingSIAM} proposed applying initial sketches (e.g., CountSketch \citep{clarkson2017LowRankApproximation} or OSNAP \citep{nelson2013OSNAP}) to each factor vector $\bm{a}_d$ in \eqref{eq:kronecker-vector}, followed by recursive pairwise combinations using TensorSketches or KFJLTs, often visualized via a binary tree structure. Extensions handle cases where $D$ is not a power of two. Song et al.~\citep{song2021FastSketchingPolynomial} adapted this recursive approach for the specific case where all factor vectors $\bm{a}_d$ are identical. Ma et al.~\citep{ma2022CostefficientGaussian} generalized the recursive structure beyond trees to arbitrary graphs, focusing on Gaussian sketches at the nodes.

\paragraph{Sampling-based sketches}
The fifth and final class comprises sampling-based methods. Unlike the predominantly \textit{oblivious} sketches mentioned above (which are independent of the data matrix $\bm{A}$), sampling-based sketches are \textit{non-oblivious}, utilizing information derived from $\bm{A}$—typically related to row norms or leverage scores—to construct the sketch $\bm{S}$. Such methods have been developed for matrices with various structures, including full Kronecker product matrices \citep{diao2019OptimalSketching, fahrbach2022SubquadraticKroneckera} and Khatri--Rao product matrices \citep{cheng2016SPALSFast, larsen2020PracticalLeverageBased, woodruff2020NearInputSparsity, malik2022MoreEfficientSampling, woodruff2022LeverageScoreSampling,PDE_inverseProblem_sketch2020, bharadwaj2023fast}. Additionally, Malik et al.~\citep{malik2021SamplingBasedMethod, malik2022MoreEfficientSampling} developed sampling techniques for matrices arising in tensor ring decomposition algorithms, whose columns involve sums of Kronecker products. Our work falls into this non-oblivious category, focusing specifically on efficient exact leverage score sampling for lower subsets of Kronecker products, a structure not explicitly addressed by previous sampling methods. A primary benefit of this type of sketch is that not all entries of the right-hand side $\bm{b}$ need to be observed, so this is particularly beneficial when forming the right-hand side is expensive, such as solving a differential equation for a given parameter set.

\section{Preliminaries}
\label{sec:preliminaries}

\subsection{Leverage Score Sampling}
\label{ssec:leverage-score}

We begin by defining the statistical leverage scores of a matrix.

\begin{definition}[Leverage Score] % Consider adding a canonical reference
    Given a matrix $\bm{A} \in \mathbb{R}^{M \times N}$, the leverage score of the $m$-th row, $\bm{A}_{m,:}$, is defined via the maximization problem:
    \begin{equation}\label{eq:ls-def}
        \tau_m(\bm{A}) \coloneqq \max_{\substack{\bm{x} \in \mathbb{R}^N \\ \bm{A}\bm{x} \neq \bm{0}}} \frac{[\bm{A} \bm{x}]_m^2}{\lVert \bm{A} \bm{x} \rVert_2^2} = \max_{\substack{\bm{x} \in \mathbb{R}^N \\ \lVert \bm{A}\bm{x} \rVert_2 = 1}} [\bm{A} \bm{x}]_m^2,
    \end{equation}
    where $[\bm{v}]_m$ denotes the $m$-th entry of vector $\bm{v}$.
\end{definition}

Leverage scores quantify the influence of each row on the column space of $\bm{A}$. If $\bm{A}$ has full column rank (i.e., $\text{rank}(\bm{A}) = N \le M$), the leverage scores sum exactly to the rank:
\begin{equation}
    \sum_{m=1}^{M} \tau_m(\bm{A}) = N.
\end{equation}
Since $\tau_m(\bm{A}) \ge 0$ by definition, normalizing the scores yields a probability distribution over the rows of $\bm{A}$:
\begin{equation}
    p_m = \mathbb{P}[I = m] = \frac{\tau_m(\bm{A})}{N}, \quad m = 1, \dots, M.
\end{equation}
This distribution $p = (p_1, \dots, p_M)$ is fundamental for constructing randomized matrix algorithms, such as row sampling or sketching operators $\bm{S}$ \citep{drineas2012FastApproximation}. Here, $p_m$ represents the probability assigned to selecting the $m$-th row.

A key property relates leverage scores to any orthonormal basis for the column space of $\bm{A}$.
\begin{proposition}[Equation 4.4 in \citep{meyer2023near}]
\label{prop:leverage-score}
    Let $\bm{Q} \in \mathbb{R}^{M \times N}$ be any matrix whose columns form an orthonormal basis for the column space of $\bm{A}$ (i.e., $\textup{span}(\bm{Q}) = \textup{span}(\bm{A})$ and $\bm{Q}^\mathsf{T} \bm{Q} = \bm{I}_N$). Then the leverage score of the $m$-th row is  the squared Euclidean norm of the $m$-th row of $\bm{Q}$:
    \begin{equation}
        \tau_m(\bm{A}) = \sum_{n=1}^{N} q_{mn}^2 = \| \bm{Q}_{m,:} \|_2^2,
    \end{equation}
    where $q_{mn}$ denotes the $(m, n)$-th entry of $\bm{Q}$ and $\bm{I}_N$ is the $N \times N$ identity matrix; and this leverage score is well-defined, meaning it is invariant under the specific $\bm{Q}$ as long as $\bm{Q}$ satisifes the above properties.
\end{proposition}

This proposition implies that leverage scores can be computed from the $\bm{Q}$ factor of the thin QR decomposition of $\bm{A}$. However, computing this decomposition requires $\mathcal{O}(MN^2)$ floating-point operations, which can be computationally prohibitive for large matrices, particularly when $M$ is large. Fortunately, intrinsic matrix structures can often be exploited to accelerate this computation, as discussed next for Kronecker product structures.

\subsection{Kronecker Product Structure and Lower Column Subsets}
\label{ssec:kronecker_structure}

Consider matrices constructed via the Kronecker product. Let $\bm{A}^{(d)} \in \mathbb{R}^{M_d \times N_d}$ for $d = 1, \dots, D$ be factor matrices, each admitting a thin QR decomposition $\bm{A}^{(d)} = \bm{Q}^{(d)} \bm{R}^{(d)}$, where $\bm{Q}^{(d)} \in \mathbb{R}^{M_d \times N_d}$ has orthonormal columns. The full Kronecker product matrix is defined as:
\begin{equation}
    \bm{A}_{\textup{kr}} = \bigotimes_{d=1}^D \bm{A}^{(d)} \in \mathbb{R}^{\left(\prod_{d=1}^D M_d\right) \times \left(\prod_{d=1}^D N_d\right)}.
\end{equation}
A well-known property is that the orthonormal basis matrix $\bm{Q}_{\textup{kr}}$ for the column space of $\bm{A}_{\textup{kr}}$ is given by the Kronecker product of the individual basis matrices:
\begin{equation}
    \bm{Q}_{\textup{kr}} = \bigotimes_{d=1}^D \bm{Q}^{(d)}.
\end{equation}
Consequently, computing the orthonormal basis $\bm{Q}_{\textup{kr}}$ leverages the structure: instead of a direct QR decomposition of $\bm{A}_{\textup{kr}}$ costing $ \mathcal{O}\left(\left(\prod_{d=1}^D M_d\right) \left(\prod_{d=1}^D N_d\right)^2\right)$, one computes the $D$ smaller QR decompositions at a significantly reduced total cost of $\mathcal{O}\left(\sum_{d=1}^D M_d N_d^2\right)$; furthermore, one does not need to even compute all $\prod_{d=1}^D M_d$ leverage scores but can instead implicitly sample from them efficiently.

In many applications, the matrix of interest, $\bm{A}$, is not the full Kronecker product $\bm{A}_{\textup{kr}}$ but rather comprises a subset of its columns. We focus on the case where $\bm{A}$ corresponds to columns indexed by a specific type of multi-index set. To the best of our knowledge, there is no prior known efficient leverage score sampling for this case.

\begin{definition}[Multi-Index Set]
    Let $\bm{N} = (N_1, \dots, N_D) \in \mathbb{N}^D$ be the tuple of column dimensions of the factor matrices $\bm{A}^{(d)}$. Define the index range for each dimension as $[N_d] \coloneqq \{1, 2, \dots, N_d\}$. The set of all possible column multi-indices for $\bm{A}_{\textup{kr}}$ is the Cartesian product
    \begin{equation}
        [\bm{N}] \coloneqq [N_1] \times \dots \times [N_D] = \prod_{d=1}^D [N_d] \subset \mathbb{N}^D.
    \end{equation}
    The total number of columns in $\bm{A}_{\textup{kr}}$ is $|[\bm{N}]| = \prod_{d=1}^D N_d$.
\end{definition}

We consider matrices $\bm{A}$ formed by selecting columns of $\bm{A}_{\textup{kr}}$ indexed by a subset $\mathcal{J} \subseteq [\bm{N}]$.
\begin{definition}[Lower Subset] % Add reference if appropriate
\label{def:lower-set}
    A subset $\mathcal{J} \subseteq [\bm{N}]$ is called a \emph{lower subset} (or downward-closed set) if, for any multi-index $\bm{\alpha} = (\alpha_1, \dots, \alpha_D) \in \mathcal{J}$, all multi-indices $\bm{\beta} = (\beta_1, \dots, \beta_D)$ satisfying $\bm{1} \leq \bm{\beta} \leq \bm{\alpha}$ (component-wise inequality) are also contained in $\mathcal{J}$:
    \begin{equation}
        \bm{\alpha} \in \mathcal{J} \text{ and } \bm{1} \leq \bm{\beta} \leq \bm{\alpha} \implies \bm{\beta} \in \mathcal{J}.
    \end{equation}
\end{definition}

Lower subsets naturally arise, for instance, as level sets $f^{-1}((-\infty, k])$ of component-wise non-decreasing functions $f: \mathbb{N}^D \rightarrow \mathbb{R}$, such as weighted $\ell_p$-norms.

The design matrix $\bm{A} \in \mathbb{R}^{M \times N}$ considered in this work consists of the columns of $\bm{A}_{\textup{kr}}$ corresponding to multi-indices in a lower subset $\mathcal{J}$, where $N = |\mathcal{J}|$. Assuming a standard (e.g., lexicographical) ordering mapping $[\bm{N}]$ to $\{1, \dots, \prod N_d\}$, we denote this relationship as
\begin{equation}
\label{eq:A-def}
    \bm{A} = (\bm{A}_{\textup{kr}})_{:, \mathcal{J}}.
\end{equation}
Our objective is to efficiently sample rows of $\bm{A}$ according to its leverage scores $\tau_m(\bm{A})$. While the leverage scores of the full Kronecker product $\bm{A}_{\textup{kr}}$ exhibit a multiplicative structure related to $\tau_{m_d}(\bm{A}^{(d)})$, this structure does not directly extend to the column subset $\bm{A}$. Nonetheless, the algorithm proposed in Section~\ref{ssec:sampling-algorithms} achieves efficient leverage score sampling for $\bm{A}$ when $\mathcal{J}$ is a lower subset. This specific structural assumption distinguishes our framework from several related works on Kronecker least squares \citep{fausett_large_1994, fausett_improved_1997, seshadri_kronecker_2017, marco_least_2019, fausett_overview_2020}.

\begin{figure}[htbp]
  \centering
  % Adjust width as needed, e.g., 0.9\textwidth or \columnwidth
  \includegraphics[trim={0 1cm 0 1cm}, clip, width=0.9\textwidth]{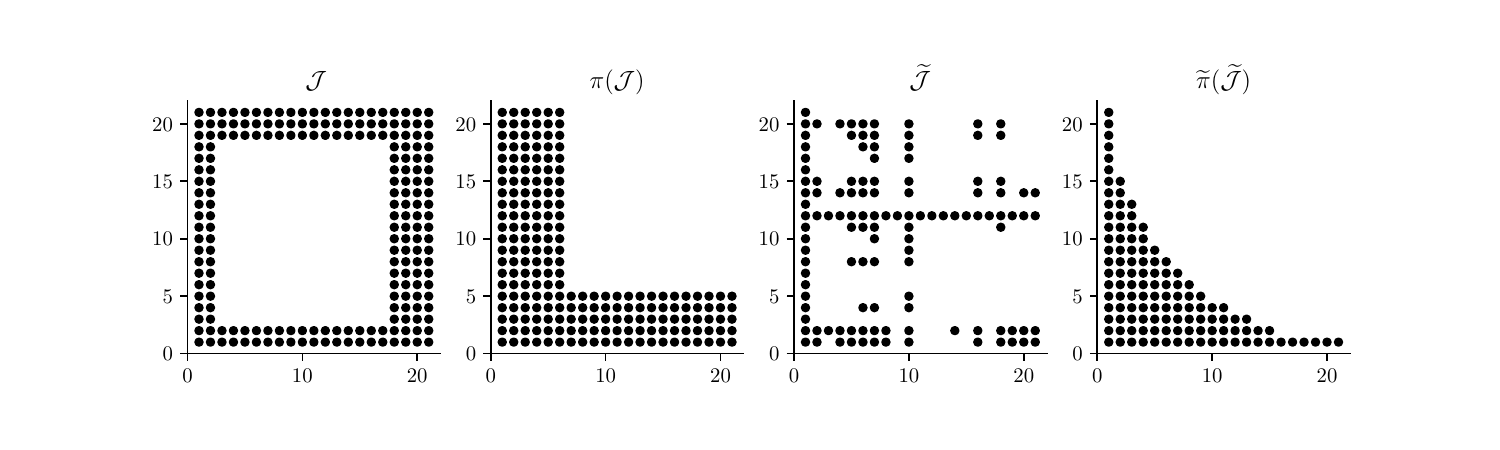}
  \caption{Illustration of index sets transformable into lower subsets. The sets $\mathcal{J}$ (left) and $\widetilde{\mathcal{J}}$ (third from left) are not lower subsets. Applying suitable dimension-wise permutations $\pi$ and $\widetilde{\pi}$ yields the sets $\pi(\mathcal{J})$ (second from left) and $\widetilde{\pi}(\widetilde{\mathcal{J}})$ (right), which are lower subsets.}
  \label{fig:permuted-sets}
\end{figure}

Furthermore, our approach extends to index sets $\mathcal{J}$ that become lower subsets after applying dimension-wise permutations. That is, if there exist permutations $\pi_d: [N_d] \rightarrow [N_d]$, for $d=1, \dots, D$, such that the transformed set $\pi(\mathcal{J}) \coloneqq \{ (\pi_1(\alpha_1), \dots, \pi_D(\alpha_D)) \mid (\alpha_1, \dots, \alpha_D) \in \mathcal{J} \}$ is a lower subset, our methods remain applicable. For instance, consider $D=2$ and the set
\begin{equation}
  \mathcal{J} = \left\{ \begin{array}{cccc} (1,1), & (1,2), & (1,3), & (1,4) \\
                                    (2,1), &        & (2,3)\hphantom{,}       &  \\
                                    (3,1), & (3,2), & (3,3), & (3,4)
                   \end{array}\right\},
\end{equation}
which is not lower. Applying permutations $\pi_1: \{1,2,3\} \mapsto \{1,3,2\}$ and $\pi_2: \{1,2,3,4\} \mapsto \{1,4,2,3\}$ yields
\begin{equation}
 \pi(\mathcal{J}) = \left\{ \begin{array}{cccc} (1,1), & (1,4), & (1,2), & (1,3) \\
                                  (3,1), &        &        & (3,2) \\
                                  (2,1), & (2,4), & (2,2), & (2,3)
                 \end{array}\right\}
                = \left\{ \begin{array}{cccc} (1,1), & (1,2), & (1,3), & (1,4) \\
                                  (2,1), & (2,2), & (2,3), & (2,4) \\
                                  (3,1), & (3,2)\hphantom{,} &        &
                 \end{array}\right\},
\end{equation}
which is a lower subset (after reordering columns for visualization). Figure~\ref{fig:permuted-sets} illustrates this concept graphically.

\section{Our Contribution: Efficient Leverage Score Sampling on Kronecker Lower Subsets}
\label{sec:effective-leverage-score}

In this section, we discuss the proposed algorithm for \emph{exact} leverage score sampling for the design matrix $\bm A$, whose columns are a lower subset of $\bm A_{\kr}$. We also present the theoretical proof showing the equivalence between the proposed algorithm and the exact leverage score sampling.

\subsection{Sampling Algorithm}
\label{ssec:sampling-algorithms}

We assume each factor matrix $\boldsymbol{A}^{(d)}$ for $d \in \{1, \dots, D\}$ has full column rank, i.e., $M_d \ge N_d$ and $\textup{rank}(\bm A^{(d)}) = N_d$. This assumption ensures that each $\boldsymbol{A}^{(d)}$ admits a thin QR decomposition:
\begin{equation}\label{eq:univariate-QR}
    \boldsymbol{A}^{(d)} = \boldsymbol{Q}^{(d)} \boldsymbol{R}^{(d)}, \quad \text{where} \quad \boldsymbol{Q}^{(d)} \in \mathbb{R}^{M_d \times N_d}, \; \boldsymbol{R}^{(d)} \in \mathbb{R}^{N_d \times N_d}.
\end{equation}
Here, $\boldsymbol{Q}^{(d)}$ has orthonormal columns, satisfying $\left(\boldsymbol{Q}^{(d)}\right)^\top \boldsymbol{Q}^{(d)} = \boldsymbol{I}_{N_d}$ (the $N_d \times N_d$ identity matrix), and $\boldsymbol{R}^{(d)}$ is upper triangular and invertible.

We assume the matrices $\boldsymbol{Q}^{(d)}$ are available, obtainable via a one-time offline preprocessing step for all $d \in \{1, \dots, D\}$. This preprocessing requires a total computational complexity of $\mathcal{O}(\sum_{d=1}^D M_{d} N_{d}^2)$ and storage complexity of $\mathcal{O}(\sum_{d=1}^D M_{d} N_{d})$. Specifically, we require access to the elements and columns of each $\boldsymbol{Q}^{(d)}$:
\begin{equation}\label{eq:Q_columns}
  \boldsymbol{Q}^{(d)} = \left[ \boldsymbol{q}_1^{(d)}, \; \ldots, \; \boldsymbol{q}_{N_d}^{(d)}\right], \quad \text{where } \boldsymbol{q}_n^{(d)} \in \mathbb{R}^{M_d}.
\end{equation}
Let $q^{(d)}_{m, n}$ denote the entry in the $m$-th row and $n$-th column of $\boldsymbol{Q}^{(d)}$, and let $\boldsymbol{Q}^{(d)}_{m,:}$ denote the $m$-th row of $\boldsymbol{Q}^{(d)}$.

The leverage scores for the matrix $\boldsymbol{A}^{(d)}$ are proportional to the squared Euclidean norms of the rows of $\boldsymbol{Q}^{(d)}$ (see Proposition~\ref{prop:leverage-score}). The probability of sampling row index $m \in \{1, \dots, M_d\}$ for the $d$-th factor matrix according to its leverage scores is given by:
\begin{equation}
\label{eq:row-prob}
    \mathbb{P}[I^{(d)} = m] = \frac{\|\boldsymbol{Q}^{(d)}_{m,:}\|_2^2}{\text{rank}(\boldsymbol{A}^{(d)})} = \frac{\|\boldsymbol{Q}^{(d)}_{m,:}\|_2^2}{N_d} = \frac{1}{N_d} \sum_{n=1}^{N_d} \left(q^{(d)}_{m, n}\right)^2,
\end{equation}
where the denominator $N_d$ arises because $\sum_{m=1}^{M_d} \|\boldsymbol{Q}^{(d)}_{m,:}\|_2^2 = \|\boldsymbol{Q}^{(d)}\|_F^2 = N_d$.

The following proposition establishes an equivalent two-stage procedure for sampling an index $I^{(d)}$ according to the distribution~\eqref{eq:row-prob}.

\begin{proposition}[Two-Step Univariate Sampling Equivalence]
\label{prop:equavalence}
    Sampling an index $I^{(d)}$ according to the probability distribution defined in Equation~\eqref{eq:row-prob} is equivalent to the following two-step process:
    \begin{enumerate}
        \item Select a column index $n$ uniformly at random from $[N_d]$.
        \item Sample a row index $m$ from $[M_d]$ according to the distribution $\mathbb{P}[I^{(d)}=m \mid J^{(d)}=n] = \left(q^{(d)}_{m, n}\right)^2$. Note that this is a valid probability distribution because $\sum_{m=1}^{M_d} (q^{(d)}_{m, n})^2 = \|\boldsymbol{q}_{n}^{(d)}\|_2^2 = 1$, since $\boldsymbol{q}_{n}^{(d)}$ is a column of the matrix $\boldsymbol{Q}^{(d)}$ with orthonormal columns.
    \end{enumerate}
\end{proposition}

The proof of Proposition~\ref{prop:equavalence} is provided in Section~\ref{ssec:proof}. % Assuming this section exists
For a single matrix $\boldsymbol{A}^{(d)}$, this proposition offers no direct advantage for online sampling efficiency, as generating samples still relies on the precomputed $\boldsymbol{Q}^{(d)}$, obtained via QR decomposition costing $\mathcal{O}(M_d N_d^2)$.

However, consider sampling from a large matrix $\boldsymbol{A}$ defined in Equation~\eqref{eq:A-def}. In this context, an extension of Proposition~\ref{prop:equavalence} to the full matrix $\boldsymbol{A}$ enables a substantial reduction in the online sampling cost compared to methods requiring explicit computation or storage of the leverage scores for $\boldsymbol{A}$. This extension motivates the main theorem of this work.

\begin{theorem}[Two-Step Multivariate Sampling Equivalence]
\label{thm:equivalence}
    Given the matrix $\boldsymbol{A}$ defined in Equation~\eqref{eq:A-def} with column multi-indices $\bm \alpha \in \Jc$, where $\Jc \subseteq [\bm N]$ satisfies the lower subset property (Definition~\ref{def:lower-set}), sampling a row multi-index $\Ibf = (I^{(1)}, \dots, I^{(D)})$ according to the leverage score distribution of $\boldsymbol{A}$ is equivalent to the following two-step procedure:
    \begin{enumerate}
        \item Sample a column multi-index $\bm \alpha = (\alpha_1, \dots, \alpha_D)$ uniformly at random from the set $\Jc$.
        \item For each dimension $d=1,\dots,D$, independently sample a row index $I^{(d)}$ from $[M_d]$ according to the probability mass function $\mathbb{P}[I^{(d)}=m \mid \alpha_{\textup{sampled}, d}=\alpha_d] = (q_{m, \alpha_d}^{(d)})^2$. The final sample is the row multi-index $\Ibf = (I^{(1)}, \dots, I^{(D)})$.
    \end{enumerate}
\end{theorem}

Theorem~\ref{thm:equivalence} establishes the validity of this efficient two-step sampling procedure for the Kronecker-structured sub-column matrix $\boldsymbol{A}$. The proof is in Section~\ref{ssec:proof}. The detailed implementation is provided in Algorithm~\ref{alg:J-lower}. This procedure significantly reduces the online sampling cost by leveraging the precomputed factors $\boldsymbol{Q}^{(d)}$ (obtained with offline complexity $\mathcal{O}(\sum_{d=1}^D M_{d} N_{d}^2)$) without needing to form or analyze the full matrix $\boldsymbol{A}$.

\begin{algorithm}[htbp]
  \KwIn{Multi-index set $\mathcal{J}$, design matrices $\bm A^{(d)}, d=1,\dots,D$.}
  \KwOut{Generated row sample $\bm I\in\mathbb{N}^D$}\Indp
  Compute $\bm Q^{(d)}$ in \eqref{eq:univariate-QR}\tcp*{$\mathcal{O}(\sum_{d=1}^DM_d N_d^2)$ cost, $\mathcal{O}(\sum_{d=1}^DM_d N_d)$ storage}
  Generate sample $\bm \alpha$ uniformly at random from column subset $\mathcal{J}$\tcp*{$\mathcal{O}(1)$ cost}
  \For{$d = 1:D$}
  {
    Generate sample $I^{(d)}$ with $\mathbb{P}[I^{(d)}=m] = \left(q^{(d)}_{m, \alpha_d}\right)^2$\tcp*{$\mathcal{O}(1)$ cost}
  }
  Generate final row sample $\bm I$ from samples $(I^{(1)}, \dots, I^{(D)})$.
  \caption{A simple and effective exact leverage score sampling algorithm when $\mathcal{J}$ is a lower set.\\Total cost: $\sum_{d=1}^D M_d N_d^2$. Total storage: $\sum_{d =1}^D M_d N_d$.}\label{alg:J-lower}
\end{algorithm}

\subsection{Proofs}
\label{ssec:proof}

\begin{proof}[Proof of Proposition~\ref{prop:equavalence}]
Let $J^{(d)}$ be an auxiliary random variable uniformly distributed over the indices $\{1, \dots, N_d\}$, representing the outcome of the first step of the procedure. Thus, its probability mass function is
\begin{equation*}
    \mathbb{P}[J^{(d)} = n] = \frac{1}{N_d}, \quad \text{for } n=1, \dots, N_d.
\end{equation*}
Let $I^{(d)}$ be the random variable representing the outcome of the second step. We define the conditional probability of $I^{(d)}=m$ given $J^{(d)}=n$, corresponding to the second step of the procedure, as
\begin{equation*}
    \mathbb{P}[I^{(d)}=m \mid J^{(d)} = n] = \frac{\left(q^{(d)}_{m, n}\right)^2}{\sum_{\ell=1}^{M_d}\left(q^{(d)}_{\ell, n}\right)^2}.
\end{equation*}
By the law of total probability, the marginal probability $\mathbb{P}[I^{(d)} = m]$ is obtained by summing over all possible values of $J^{(d)}$:
\begin{equation}
\begin{aligned}
    \mathbb{P}_\textup{TwoStep}[I^{(d)} = m] &= \sum_{n=1}^{N_d} \mathbb{P}[I^{(d)}=m \mid J^{(d)} = n] \mathbb{P}[J^{(d)} = n] \\
    &= \sum_{n=1}^{N_d} \left( \frac{\left(q^{(d)}_{m, n}\right)^2}{\sum_{\ell=1}^{M_d}\left(q^{(d)}_{\ell, n}\right)^2} \right) \left( \frac{1}{N_d} \right) \\
    &= \frac{1}{N_d} \sum_{n=1}^{N_d} \frac{\left(q^{(d)}_{m, n}\right)^2}{\sum_{\ell=1}^{M_d}\left(q^{(d)}_{\ell, n}\right)^2} = \frac{1}{N_d} \sum_{n=1}^{N_d}\left(q^{(d)}_{m, n}\right)^2.\label{eq:derived_prob_im}
\end{aligned}
\end{equation}
This derived marginal probability distribution~\eqref{eq:derived_prob_im} matches the leverage score sampling distribution, thus completing the proof.
\end{proof}

To prove the main sampling equivalence result, Theorem~\ref{thm:equivalence}, we first establish in Proposition~\ref{prop:Q-lower} that the matrix $\bm{Q}$ (constructed using columns $\bm{q}_{\bm{\alpha}}$ for $\bm{\alpha} \in \mathcal{J}$) forms an orthonormal basis for the column space of $\bm{A}$. That is, we show $\bm{Q}^\top\bm{Q} = \bm{I}_{N}$ and $\operatorname{span}(\bm{Q}) = \operatorname{span}(\bm{A})$. Then, by adapting the law of total probability argument used for the univariate case (Proposition~\ref{prop:equavalence}), we demonstrate that the output distribution of the two-step procedure detailed in Algorithm~\ref{alg:J-lower} matches the exact leverage score distribution of $\bm{A}$.

Recall the following definitions and relationships:
% a column of the Kronecker product matrix A corresponding to multi-index alpha:
\begin{equation} \label{eq:aj-def} 
\bm{a}_{\bm{\alpha}} = \bigotimes_{d=1}^D \bm{a}^{(d)}_{\alpha_d} 
\end{equation}
% the relationship between columns of A^{(d)} and Q^{(d)} via R^{(d)}:
\begin{equation} \label{eq:q1d-def} 
\bm{a}^{(d)}_{k} = \sum_{i=1}^{N_d} R^{(d)}_{i,k} \bm{q}^{(d)}_i 
\end{equation}
% the definition of the basis vectors q_alpha for the column space of A:
\begin{equation} \label{eq:q-Jlower} 
\bm{q}_{\bm{\alpha}} = \bigotimes_{d=1}^D \bm{q}^{(d)}_{\alpha_d} 
\end{equation} 

\begin{proposition} \label{prop:Q-lower}
 Assume $\mathcal{J} = \{\bm{\alpha}_1, \ldots, \bm{\alpha}_N\}$ is a lower subset (Definition~\ref{def:lower-set}), where $N = |\mathcal{J}|$. Let $\bm{Q}$ be the matrix whose columns are $\{\bm{q}_{\bm{\alpha}}\}_{\bm{\alpha} \in \mathcal{J}}$, with $\bm{q}_{\bm{\alpha}}$ defined as in Equation~\eqref{eq:q-Jlower} using the orthonormal columns $\bm{q}^{(d)}_k$ of $\bm{Q}^{(d)}$ from the QR decomposition $\bm{A}^{(d)} = \bm{Q}^{(d)} \bm{R}^{(d)}$. Then the columns of $\bm{Q}$ form an orthonormal basis for the column space $V = \operatorname{span}(\bm{A})$.
\end{proposition}

\begin{proof}[Proof of Proposition~\ref{prop:Q-lower}]
Let $\bm{Q}$ be the matrix with columns $\bm{q}_{\bm{\alpha}}$ for $\bm{\alpha} \in \mathcal{J}$. We need to show two properties:
\begin{enumerate}
    \item The columns $\{\bm{q}_{\bm{\alpha}}\}_{\bm{\alpha} \in \mathcal{J}}$ are orthonormal.
    \item The columns span the same space as the columns of $\bm{A}$, i.e., $\operatorname{span}(\bm{Q}) = \operatorname{span}(\bm{A}) = V$.
\end{enumerate}

% we prove orthonormality. Consider the inner product between two columns $\bm{q}_{\bm{\alpha}_i}$ and $\bm{q}_{\bm{\alpha}_j}$ of $\bm{Q}$, where $\bm{\alpha}_i, \bm{\alpha}_j \in \mathcal{J}$. Using the definition~\eqref{eq:q-Jlower} and the mixed-product property of inner products and Kronecker products:
% \begin{align*}
%     \langle \bm{q}_{\bm{\alpha}_i} , \bm{q}_{\bm{\alpha}_j} \rangle &= \langle \bigotimes_{d=1}^D \bm{q}^{(d)}_{[\bm{\alpha}_i]_d} , \bigotimes_{d=1}^D \bm{q}^{(d)}_{[\bm{\alpha}_j]_d} \rangle \\
%     &= \prod_{d=1}^D \langle \bm{q}^{(d)}_{[\bm{\alpha}_i]_d} , \bm{q}^{(d)}_{[\bm{\alpha}_j]_d} \rangle.
% \end{align*}
% Since the columns $\{\bm{q}^{(d)}_k\}_k$ of each $\bm{Q}^{(d)}$ are orthonormal, we have $\langle \bm{q}^{(d)}_k , \bm{q}^{(d)}_l \rangle = \delta_{k,l}$ (Kronecker delta). Therefore,
% \begin{equation*}
%     \langle \bm{q}_{\bm{\alpha}_i} , \bm{q}_{\bm{\alpha}_j} \rangle = \prod_{d=1}^D \delta_{[\bm{\alpha}_i]_d, [\bm{\alpha}_j]_d}.
% \end{equation*}
% The product of Kronecker deltas is 1 if and only if $[\bm{\alpha}_i]_d = [\bm{\alpha}_j]_d$ for all $d=1, \dots, D$, which means $\bm{\alpha}_i = \bm{\alpha}_j$. Otherwise, the product is 0. Thus,
% \begin{equation*}
%     \langle \bm{q}_{\bm{\alpha}_i} , \bm{q}_{\bm{\alpha}_j} \rangle = \delta_{\bm{\alpha}_i, \bm{\alpha}_j} = \begin{cases} 1 & \text{if } i = j \\ 0 & \text{if } i \neq j \end{cases}.
% \end{equation*}
% This confirms that the $N$ columns of $\bm{Q}$ are orthonormal, satisfying $\bm{Q}^\top \bm{Q} = \bm{I}_N$.

First, the orthonormality of $\bm Q$ directly stems from the fact that $\bm Q = (\bm Q_\kr)_{:,\Jc}$ from Equation~\eqref{eq:A-def}, so we only need to show that $\operatorname{span}(\bm{Q}) = \operatorname{span}(\bm{A})$.
Since $\{\bm{q}_{\bm{\alpha}}\}_{\bm{\alpha} \in \mathcal{J}}$ is a set of $N$ orthonormal vectors and $\dim(V) = \dim(\operatorname{span}(\bm{A})) = N = |\mathcal{J}|$, it suffices to show that every column of $\bm{A}$ lies within the span of the columns of $\bm{Q}$.
Let $\bm{a}_{\bm{\alpha}}$ be an arbitrary column of $\bm{A}$ corresponding to the multi-index $\bm{\alpha} \in \mathcal{J}$. From its definition~\eqref{eq:aj-def} and the QR decomposition relationship~\eqref{eq:q1d-def}, we have:
\begin{align*}
    \bm{a}_{\bm{\alpha}} &= \bigotimes_{d=1}^D \bm{a}^{(d)}_{\alpha_d} \\
    &= \bigotimes_{d=1}^D \left( \sum_{i_d=1}^{N_d} R^{(d)}_{i_d, \alpha_d} \bm{q}^{(d)}_{i_d} \right) \\
    &= \sum_{i_1=1}^{N_1} \dots \sum_{i_D=1}^{N_D} \left( \prod_{d=1}^D R^{(d)}_{i_d, \alpha_d} \right) \left( \bigotimes_{d=1}^D \bm{q}^{(d)}_{i_d} \right) \quad \text{(by multi-linearity of } \otimes \text{)} \\
    &= \sum_{\bm{\beta} \in [\bm N]} \left( \prod_{d=1}^D R^{(d)}_{\beta_d, \alpha_d} \right) \bm{q}_{\bm{\beta}}.
\end{align*}
Since each $\bm{R}^{(d)}$ is upper triangular, the coefficient $R^{(d)}_{\beta_d, \alpha_d}$ is zero if $\beta_d > \alpha_d$. Thus, the product $\prod_{d=1}^D R^{(d)}_{\beta_d, \alpha_d}$ is non-zero only if $\beta_d \le \alpha_d$ for all $d$ (denoted $\bm{\beta} \le \bm{\alpha}$). The summation is therefore restricted:
\begin{equation*}
    \bm{a}_{\bm{\alpha}} = \sum_{\bm{\beta} \le \bm{\alpha}} \left( \prod_{d=1}^D R^{(d)}_{\beta_d, \alpha_d} \right) \bm{q}_{\bm{\beta}}.
\end{equation*}
By the definition of a lower set (Definition~\ref{def:lower-set}), if $\bm{\alpha} \in \mathcal{J}$ and $\bm{\beta} \le \bm{\alpha}$, then $\bm{\beta} \in \mathcal{J}$. Thus, the sum is effectively over indices $\bm{\beta}$ in $\mathcal{J}$:
\begin{equation*}
    \bm{a}_{\bm{\alpha}} = \sum_{\substack{\bm{\beta} \in \mathcal{J} \\ \bm{\beta} \le \bm{\alpha}}} \left( \prod_{d=1}^D R^{(d)}_{\beta_d, \alpha_d} \right) \bm{q}_{\bm{\beta}} = \sum_{\bm{\beta} \in \mathcal{J}} c_{\bm{\alpha}, \bm{\beta}} \, \bm{q}_{\bm{\beta}},
\end{equation*}
where $c_{\bm{\alpha}, \bm{\beta}} = \prod_{d=1}^D R^{(d)}_{\beta_d, \alpha_d}$ if $\bm{\beta} \le \bm{\alpha}$, and $c_{\bm{\alpha}, \bm{\beta}} = 0$ otherwise.

This shows that every column $\bm{a}_{\bm{\alpha}}$ of $\bm{A}$ is a linear combination of the columns $\{\bm{q}_{\bm{\beta}}\}_{\bm{\beta} \in \mathcal{J}}$ of $\bm{Q}$. Therefore, $\operatorname{span}(\bm{A}) \subseteq \operatorname{span}(\bm{Q})$. Since both spaces have the same dimension $N$, we conclude that $\operatorname{span}(\bm{A}) = \operatorname{span}(\bm{Q})$.

Combining the orthonormality ($\bm{Q}^\top\bm{Q} = \bm{I}_N$) and the spanning property, the set of columns $\{\bm{q}_{\bm{\alpha}}\}_{\bm{\alpha} \in \mathcal{J}}$ forms an orthonormal basis for $V = \operatorname{span}(\bm{A})$.
\end{proof}

% --- Proof for Theorem 1 ---
\begin{proof}[Proof of Theorem~\ref{thm:equivalence}]
The goal is to show that the probability distribution generated by the two-step sampling procedure in Algorithm~\ref{alg:J-lower} matches the leverage score sampling distribution for the full matrix $\bm{A}$. Let $N = |\mathcal{J}| = \operatorname{rank}(\bm{A})$.

First, consider the leverage score distribution. The leverage score of the row corresponding to the multi-index $\bm{m} = (m_1, \dots, m_D)$ of $\bm{A}$ is proportional to the squared Euclidean norm of the $\bm{m}$-th row of the matrix $\bm{Q}$ whose columns form an orthonormal basis for the column space of $\bm{A}$. As established in Proposition~\ref{prop:Q-lower}, the columns of $\bm{Q}$ are $\{\bm{q}_{\bm{\alpha}}\}_{\bm{\alpha} \in \mathcal{J}}$, where $\bm{q}_{\bm{\alpha}} = \bigotimes_{d=1}^D \bm{q}^{(d)}_{\alpha_d}$. The $\bm{m}$-th entry of the column vector $\bm{q}_{\bm{\alpha}}$ is given by the product of the corresponding entries from the factor vectors: $(\bm{q}_{\bm{\alpha}})_{\bm{m}} = \prod_{d=1}^D q^{(d)}_{m_d, \alpha_d}$.

The squared norm of the $\bm{m}$-th row of $\bm{Q}$, denoted $\bm{Q}_{\bm{m}, :}$, is the sum of the squares of its entries:
\begin{equation*}
    \|\bm{Q}_{\bm{m}, :}\|_2^2 = \sum_{\bm{\alpha} \in \mathcal{J}} \left( (\bm{q}_{\bm{\alpha}})_{\bm{m}} \right)^2 = \sum_{\bm{\alpha} \in \mathcal{J}} \left( \prod_{d=1}^D q^{(d)}_{m_d, \alpha_d} \right)^2.
\end{equation*}
The probability of sampling the row multi-index $\bm{m}$ according to the leverage scores is this squared norm divided by the rank $N$:
\begin{equation} \label{eq:lev_score_prob_multi}
    \mathbb{P}_{\text{LevScore}}[\bm{I} = \bm{m}] = \frac{\|\bm{Q}_{\bm{m}, :}\|_2^2}{N} = \frac{1}{N} \sum_{\bm{\alpha} \in \mathcal{J}} \prod_{d=1}^D \left( q^{(d)}_{m_d, \alpha_d} \right)^2.
\end{equation}

Now, consider the two-step sampling procedure described in Theorem~\ref{thm:equivalence}.
\begin{enumerate}
    \item Let $\bm{\alpha}_{\text{sampled}} = (\alpha_{\text{sampled},1}, \dots, \alpha_{\text{sampled},D})$ be the random multi-index sampled uniformly from $\mathcal{J}$ in Step 1. Its probability mass function is $\mathbb{P}[\bm{\alpha}_{\text{sampled}} = \bm{\alpha}] = 1/N$ for $\bm{\alpha} \in \mathcal{J}$.
    \item Let $\bm{I} = (I^{(1)}, \dots, I^{(D)})$ be the random multi-index generated in Step 2. Given $\bm{\alpha}_{\text{sampled}} = \bm{\alpha}$, the components $I^{(d)}$ are sampled independently for $d=1, \dots, D$, with $\mathbb{P}[I^{(d)} = m_d \mid \alpha_{\text{sampled},d} = \alpha_d] = (q^{(d)}_{m_d, \alpha_d})^2$. The conditional probability of obtaining the specific multi-index $\bm{m} = (m_1, \dots, m_D)$ given $\bm{\alpha}_{\text{sampled}} = \bm{\alpha}$ is:
    \begin{equation*}
        \mathbb{P}[\bm{I} = \bm{m} \mid \bm{\alpha}_{\text{sampled}} = \bm{\alpha}] = \prod_{d=1}^D \mathbb{P}[I^{(d)} = m_d \mid \alpha_{\text{sampled},d} = \alpha_d] = \prod_{d=1}^D \left(q^{(d)}_{m_d, \alpha_d}\right)^2.
    \end{equation*}
    Note that the normalization $\sum_{m_d=1}^{M_d} \left(q^{(d)}_{m_d, \alpha_d}\right)^2 = \|\bm{q}^{(d)}_{\alpha_d}\|_2^2 = 1$ ensures this is a valid conditional probability for each $d$.
\end{enumerate}
Using the law of total probability, we find the marginal probability distribution for $\bm{I}$ generated by this two-step procedure:
\begin{align} 
    \mathbb{P}_{\text{TwoStep}}[\bm{I} = \bm{m}] &= \sum_{\bm{\alpha} \in \mathcal{J}} \mathbb{P}[\bm{I} = \bm{m} \mid \bm{\alpha}_{\text{sampled}} = \bm{\alpha}] \mathbb{P}[\bm{\alpha}_{\text{sampled}} = \bm{\alpha}] \nonumber \\
    &= \sum_{\bm{\alpha} \in \mathcal{J}} \left( \prod_{d=1}^D (q^{(d)}_{m_d, \alpha_d})^2 \right) \left( \frac{1}{N} \right) \nonumber \\
    &= \frac{1}{N} \sum_{\bm{\alpha} \in \mathcal{J}} \prod_{d=1}^D (q^{(d)}_{m_d, \alpha_d})^2. \label{eq:two_step_prob_multi_final} 
\end{align}
Comparing the resulting distribution~\eqref{eq:two_step_prob_multi_final} with the leverage score distribution~\eqref{eq:lev_score_prob_multi}, we see they are identical: $\mathbb{P}_{\text{TwoStep}}[\bm{I} = \bm{m}] = \mathbb{P}_{\text{LevScore}}[\bm{I} = \bm{m}]$. Therefore, the two-step sampling procedure is equivalent to sampling according to the leverage scores of $\bm{A}$.
\end{proof}

\section{Experiments}
\label{sec:experiments}

In this section, we numerically demonstrate that the exact leverage score sampling method proposed in Algorithm~\ref{alg:J-lower} yields improved accuracy compared to uniform random row sampling and an approximate leverage score sampling approach based on the full tensor product matrix. The test cases are drawn from polynomial approximation, specifically parametric uncertainty quantification using polynomial chaos expansions (PCE) \citep{ghanem2003stochastic}. We utilize point-weight pairs from $M_d$-point Gauss-Legendre quadrature on $[-1,1]$ and employ one-dimensional Legendre polynomials $\{a_n^{(d)}\}_{n \in [N_d]}$ as basis functions within each dimension $d$.

For a given multi-index set $\mathcal{J}$, the least squares matrix $\bm{A}$ is constructed as the $\mathcal{J}$-column subset of the full tensor product matrix $\bm{A}_{\textup{kr}}$, resulting in $M = \prod_{d=1}^D M_d$ rows. We consider two standard choices for the index set $\mathcal{J}$: total degree (TD, $\mathcal{J} = \mathcal{J}_\textup{TD}$) and hyperbolic cross (HC, $\mathcal{J} = \mathcal{J}_\hc$) sets with space degree $k$, both of which are lower subsets as defined in Definition~\ref{def:lower-set}. In our examples, the quadrature nodes, weights, and basis functions are identical across all dimensions $d \in \{1, \dots, D\}$, implying that the factor matrices $\bm{A}^{(1)}, \ldots, \bm{A}^{(D)}$ are identical. The specific definition of the right-hand side vector $\bm{b}$ will be provided for each example.

We perform sketched least squares using Algorithm~\ref{alg:J-lower} and compare the performance of the following three row sampling strategies for constructing the sketch:
\begin{itemize}
  \item \textbf{Uniform Sampling:} Rows are selected by sampling row multi-indices $(m_1, \dots, m_D)$ uniformly at random, where each $m_d$ is chosen uniformly from $\{1, \dots, M_d\}$. Generating each random row multi-index requires $\mathcal{O}(D)$ operations. This strategy does not utilize any structural information from $\bm{A}$.
  \item \textbf{TP Leverage Sampling (Approximate):} Rows are sampled according to the leverage scores of the full tensor product matrix $\bm{A}_{\textup{kr}} = \bigotimes_{d=1}^D \bm{A}^{(d)}$. This serves as an easily computable approximation to the true leverage scores of $\bm{A}$, implicitly assuming the scores of $\bm{A}_{\textup{kr}}$ are representative (e.g., when $\mathcal{J}$ is close to the full index set $[\bm{N}]$). This method requires the QR factors of the $\bm{A}^{(d)}$ matrices, computable with a one-time setup cost of $\mathcal{O}(\sum_{d=1}^D M_d N_d^2)$. Subsequent sampling based on these factors is computationally efficient.
  \item \textbf{Exact Leverage Sampling (Proposed):} Rows are sampled according to the exact leverage scores of the subset matrix $\bm{A}$ using our proposed Algorithm~\ref{alg:J-lower}. This method also requires the QR factors of the $\bm{A}^{(d)}$ matrices (identical $\mathcal{O}(\sum_{d=1}^D M_d N_d^2)$ setup cost as TP sampling) but employs the correct sampling distribution derived specifically for the lower subset $\mathcal{J}$. The sampling process itself is designed to be efficient after the initial setup.
\end{itemize}

For each sampling method, we compute the relative residual error of the sketched least squares solution $\tilde{\bm{x}}$, defined as
\begin{equation}
    \text{Relative error} = \frac{\| \bm{A} \tilde{\bm{x}} - \bm{b} \|_2}{\|\bs{b}\|_2}.
\end{equation}
We compare this against the optimal relative error achieved by the full least squares solution $\bm{x}^\ast$, given by
\begin{equation}
    \text{Optimal relative error} = \frac{\| \bm{A} \bm{x}^\ast - \bm{b} \|_2}{\|\bs{b}\|_2}.
\end{equation}
Since the sketched solution $\tilde{\bm{x}}$ depends on the random sketch $\bm{S}$, its relative error is a random variable. We report its empirical distribution based on 100 independent trials for each experimental setup.

We consider three distinct examples to evaluate the methods:
\begin{enumerate}
    \item A Duffing oscillator under free vibration, representing a low-dimensional problem ($D=3$) where higher-degree polynomials may be used.
    \item The Ishigami function, another benchmark problem in UQ ($D=3$).
    \item Prediction of the remaining useful life (RUL) of a Lithium-ion battery, representing a higher-dimensional problem ($D=7$) typically tackled with lower-degree polynomials.
\end{enumerate}
The code used to generate the results for these experiments is publicly available at \href{https://github.com/CU-UQ/monotone-lower-set}{https://github.com/CU-UQ/monotone-lower-set}.

\subsection{Nonlinear Duffing oscillator}
\label{ssec:empirical-duffing}

Our first example considers modeling the uncertainty in the displacement solution $u(\ybf,t)$ for a nonlinear single-degree-of-freedom Duffing oscillator under free vibration, following \citep{MaiChuVPCEF}. The system dynamics are described by:
\begin{equation}
\begin{aligned}
& \ddot{u}(\ybf, t) + 2 \omega_1 \omega_2 \dot{u}(\ybf, t) + \omega_1^2 ( u(\ybf, t) + \omega_3 u^3(\ybf, t) )= 0, \\
& {u}(\ybf, 0) = 1, \qquad \dot{u}(\ybf, 0) = 0,
\end{aligned}
\label{eqn:duffing}
\end{equation}
where the parameters $\{ \omega_d \}_{d=1}^3$ depend on three independent uncertain inputs $\{ y^{(d)} \}_{d=1}^3$, each assumed uniformly distributed in $[-1,1]$:
\begin{equation}
\label{eqn:duffing_rvs}
\omega_1 = 2 \pi (1 + 0.2 y^{(1)}), \quad
 \omega_2 = 0.05 (1+0.05 y^{(2)}), \quad
 \omega_3 = -0.5 (1+0.5 y^{(3)}).
\end{equation}
This constitutes a $D=3$ dimensional uncertainty quantification problem. We employ $M_d = 20$ Gauss-Legendre quadrature points per dimension for numerical integration based on tensor products. The quantity of interest (QoI) for the PCE model is the displacement at time $t=4$, i.e., $u(\ybf,4)$.

We construct PCE approximations using different polynomial spaces, defined by their index sets $\mathcal{J}$. Specifically, we consider TD spaces of polynomial order $k=9$ and $k=12$ (denoted $\Jc_\textup{TD}(9)$ and $\Jc_\textup{TD}(12)$) and HC spaces of order $k=15$ and $k=18$ (denoted $\Jc_\hc(15)$ and $\Jc_\hc(18)$). Higher orders are chosen for the HC spaces relative to the TD spaces to achieve roughly comparable basis set sizes ($|\mathcal{J}|$) for the different space types. For each space $\mathcal{J}$, the sketch size (number of rows sampled, $K$) is set to four times the number of basis functions, i.e., $K = 4 |\mathcal{J}|$.

\begin{figure}[ht]
	\centering
	\includegraphics[width=1.0\textwidth]{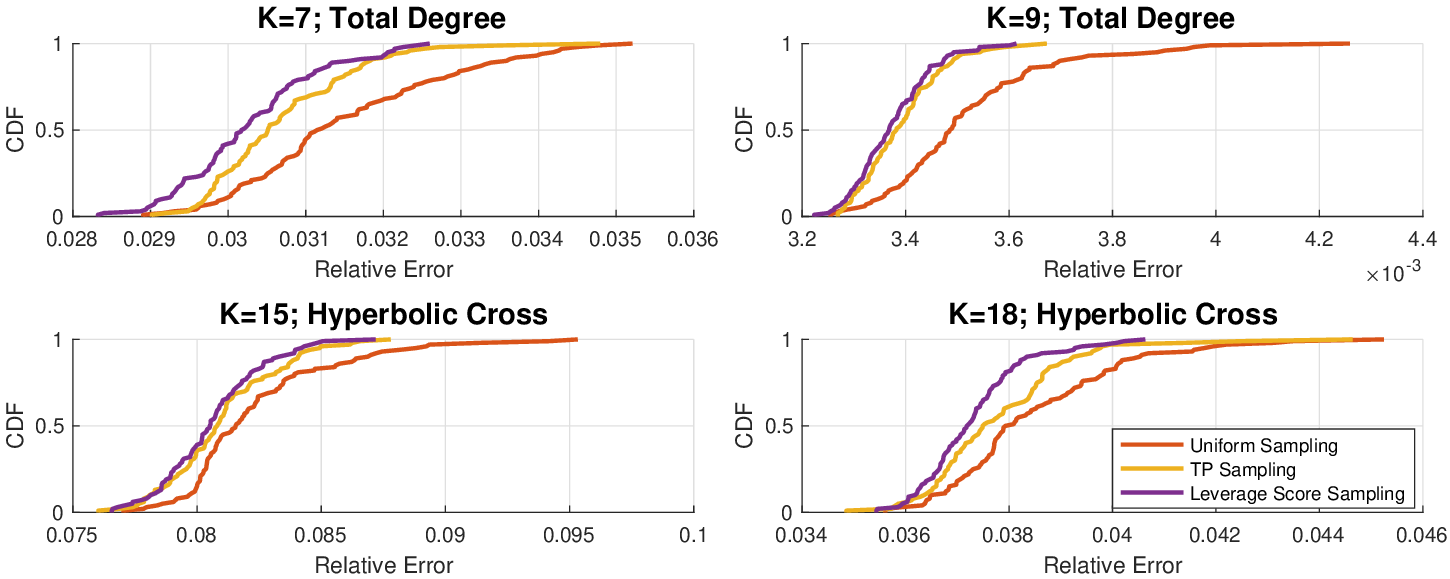}
	\caption{Empirical CDFs of the relative error in estimating the displacement $u(\ybf,4)$ of the Duffing oscillator. Top row: Total degree (TD) PCE of order $k=9$ (left, optimal rel. error $2.6\times 10^{-2}$) and $k=12$ (right, optimal rel. error $2.9\times 10^{-3}$). Bottom row: Hyperbolic cross (HC) PCE of order $k=15$ (left, optimal rel. error $6.9\times 10^{-2}$) and $k=18$ (right, optimal rel. error $3.2\times 10^{-2}$). Results compare Uniform (blue), TP Leverage (orange), and Exact Leverage (green) sampling over 100 trials.}
	\label{fig:duffing}
\end{figure}

Figure~\ref{fig:duffing} presents the empirical cumulative distribution functions (CDFs) of the relative errors obtained over 100 independent trials for the three sampling strategies (Uniform sampling, TP sampling, Leverage score sampling) across the four selected polynomial spaces. The results clearly show that exact leverage score sampling consistently outperforms uniform sampling, yielding significantly lower errors. Furthermore, exact leverage score sampling achieves higher accuracy (lower relative errors) compared to the approximate TP leverage sampling method. This aligns with expectations, as the proposed method uses the true leverage score distribution for the subset matrix $\bm{A}$, whereas TP sampling uses an approximation based on the full tensor product matrix $\bm{A}_{\textup{kr}}$. These empirical findings demonstrate the practical benefit of using exact leverage score sampling for improved accuracy in sketched least squares solutions for this class of problems.

\subsection{Ishigami function}
\label{ssec:empirical-ishigami}

The second example utilizes the Ishigami function, a common benchmark in uncertainty quantification studies \citep{IshigamiT1990Aiqt,Diaz2018}, defined as:
\begin{equation} \label{eq:ishigami}
   f(\ybf)=\sin (\pi y^{(1)})+a\sin^2 (\pi y^{(2)})+b(\pi y^{(3)})^4\sin (\pi y^{(1)}),
\end{equation}
with constants chosen as $a=7$ and $b=0.1$. The inputs $\{y^{(d)}\}_{d=1}^3$ are independent variables, each assumed uniformly distributed in $[-1,1]$, making this a $D=3$ dimensional problem. As in the previous example, we use $M_d=20$ Gauss-Legendre points per dimension for constructing the least squares system via tensor-product quadrature.

We approximate $f(\ybf)$ using PCE models based on TD polynomial spaces of order $k=7$ ($|\Jc_\textup{TD}(7)|=120$ basis functions) and $k=9$ ($|\Jc_\textup{TD}(9)|=220$), and HC spaces of order $k=15$ ($|\Jc_\hc(15)|=110$) and $k=18$ ($|\Jc_\hc(18)|=134$). Consistent with the setup in Section~\ref{ssec:empirical-duffing}, the sketch size $K$ (number of sampled rows) for the least squares solve is set to four times the basis size for each space, i.e., $K = 4 |\mathcal{J}|$.

\begin{figure}[ht]
\centering
\includegraphics[width=1.0\textwidth]{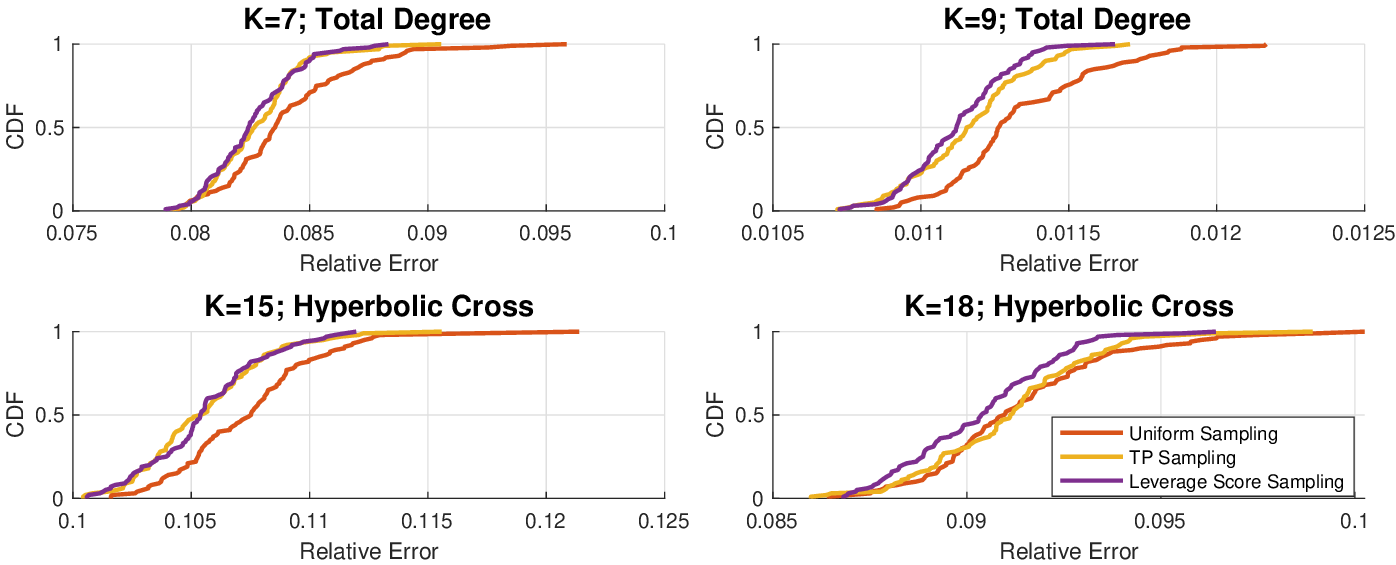}
\caption{Empirical CDFs of the relative error in approximating the Ishigami function $f(\ybf)$. Top row: Total degree (TD) PCE of order $k=7$ (left, optimal rel. error $7.0\times 10^{-3}$) and $k=9$ (right, optimal rel. error $9.5\times 10^{-4}$). Bottom row: Hyperbolic cross (HC) PCE of order $k=15$ (left, optimal rel. error $9.0\times 10^{-2}$) and $k=18$ (right, optimal rel. error $7.7\times 10^{-2}$). Results compare Uniform (blue), TP Leverage (orange), and Exact Leverage (green) sampling over 100 trials.}
\label{fig:ishigami}
\end{figure}

Figure~\ref{fig:ishigami} displays the empirical CDFs of the relative errors over 100 independent trials for the three sampling methods. The results are qualitatively consistent with those observed for the Duffing oscillator (Section~\ref{ssec:empirical-duffing}). As expected, exact leverage score sampling yields the highest accuracy among the tested methods, outperforming the approximate TP leverage sampling approach, which in turn outperforms simple uniform sampling. This further corroborates the benefit of using exact leverage scores when feasible.

\subsection{Prediction of remaining useful life of batteries}
\label{ssec:empirical-battery}

Our final example addresses a higher-dimensional problem: model-based estimation of the remaining useful life (RUL) of a Lithium-ion battery (LIB), drawing upon models described in \citep{Sankararaman2013UncertaintyQI, Sankararaman2014UQfirstOrder, Sankararaman2015significance, hadigol2018least}. RUL prediction estimates the time until a battery's health reaches a predefined threshold, often defined by its capacity degradation. The system is modeled using a state-space representation:
\begin{equation}
\begin{aligned}
	\dot{\bm{z}}(t) &= \bm{f}(t, \bm{z}(t), \bm{\theta}(t), \bm{\nu}(t), \bm{v}_p(t)), \\ % State equation
	\bm{w}(t) &= \bm{h}(t, \bm{z}(t), \bm{\theta}(t), \bm{\nu}(t), \bm{v}_m(t)), % Output equation
\end{aligned}
\end{equation}
where $\bm{z}(t)$ denotes the state vector, $\bm{\theta}(t)$ the model parameters, $\bm{\nu}(t)$ the inputs, $\bm{w}(t)$ the measured outputs, $\bm{v}_p(t)$ process noise, $\bm{v}_m(t)$ measurement noise, and $\bm{f}, \bm{h}$ are the state and output equations, respectively. For detailed model equations and parameters specific to the LIB application, we refer the interested reader to the cited works.

In this specific application, we consider $D=7$ uncertain input parameters: the constant discharge current (modeled as uniform on a given range), initial estimates for three key state variables, and three associated process noise terms. These physical inputs are transformed linearly into the canonical interval $[-1,1]$ for compatibility with the Legendre polynomial basis. Due to the higher dimensionality ($D=7$), we employ a lower polynomial degree, $k=3$, for constructing PCE approximations using both TD and HC index sets ($\Jc_\textup{TD}(3)$, $\Jc_\hc(3)$). We investigate the performance using either $M_d=4$ or $M_d=5$ Gauss-Legendre points per dimension for the tensor-product quadrature underlying the least squares problem. The sketch size (number of rows sampled) is again set relative to the number of basis functions $|\mathcal{J}|$ as $K = 4 |\mathcal{J}|$.

\begin{figure}[ht]
	\centering
	\includegraphics[width=1.0\textwidth]{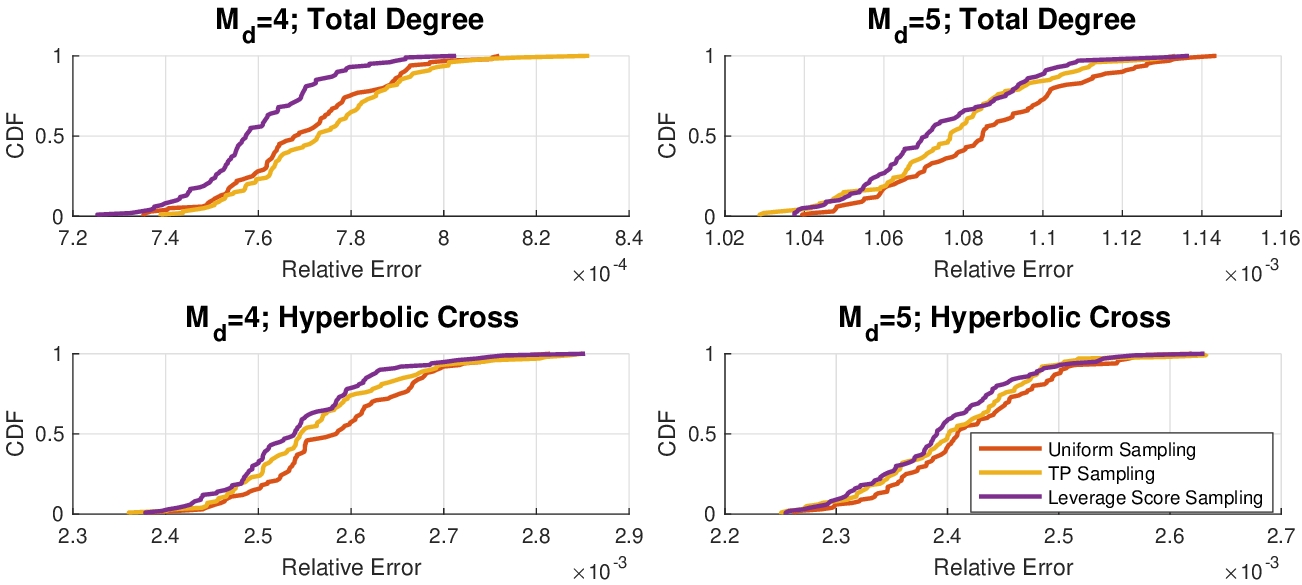}
	\caption{Empirical CDFs of the relative error in estimating the battery RUL using PCE of degree $k=3$. Top row: Total degree (TD) space with $M_d=4$ quadrature points per dimension (left, optimal rel. error $6.6\times 10^{-4}$) and $M_d=5$ (right, optimal rel. error $9.3\times 10^{-4}$). Bottom row: Hyperbolic cross (HC) space with $M_d=4$ (left, optimal rel. error $2.2\times 10^{-3}$) and $M_d=5$ (right, optimal rel. error $1.1\times 10^{-3}$). Results compare Uniform (blue), TP Leverage (orange), and Exact Leverage (green) sampling over 100 trials.}
	\label{fig:battery_0}
\end{figure}

Figure~\ref{fig:battery_0} presents the empirical CDFs of the relative errors over 100 independent trials. A notable observation occurs in the top left panel (TD space, $M_d=4$): the performance of uniform sampling is nearly identical to that of the approximate TP leverage sampling. This phenomenon arises because when the number of quadrature points $M_d$ equals the number of basis functions required for degree $k$ in one dimension ($N_d = k+1 = 4$), the corresponding factor matrix $\bm{A}^{(d)}$ becomes square. 
% For Legendre polynomials evaluated at Gauss-Legendre points, the properties of Gaussian quadrature imply that 
Therefore, the matrix $\bm{Q}^{(d)}$ in the QR decomposition of the (appropriately weighted) $\bm{A}^{(d)}$ is exactly square and orthonormal, and the leverage scores of the full tensor product matrix $\bm{A}_{\textup{kr}} = \bigotimes \bm{A}^{(d)}$ become uniform \citep{hampton_coherence_2015}, causing the TP leverage sampling strategy to degenerate to uniform sampling in this specific configuration ($M_d=k+1$). Our proposed exact leverage score sampling algorithm (Algorithm~\ref{alg:J-lower}) is not subject to this degeneracy, as it computes the (generally non-uniform) leverage scores corresponding to the actual column subset matrix $\bm{A}$ defined by $\mathcal{J}$.

Aside from this specific scenario where $M_d=k+1$, the results for the RUL problem are consistent with those observed in the previous examples (Sections~\ref{ssec:empirical-duffing} and \ref{ssec:empirical-ishigami}): exact leverage score sampling provides the best accuracy, followed by TP leverage sampling, with uniform sampling performing worst. The advantage of exact leverage scores is particularly evident for the hyperbolic cross spaces and when $M_d > k+1$.

\section{Conclusion}
\label{sec:conclusion}

Leverage score sampling offers a powerful approach for sketching large least squares problems, but the cost of using exact scores is often prohibitive. This paper introduced an efficient algorithm (Algorithm~\ref{alg:J-lower}) for exact leverage score sampling, specifically tailored to matrices that form lower column subsets of Kronecker products—a structure prevalent in high-dimensional approximation and uncertainty quantification. By exploiting this inherent structure, our method bypasses the need for expensive full-matrix leverage score evaluations and makes exact sampling computationally feasible for this class.

Numerical experiments involving polynomial chaos expansion models (Duffing oscillator, Ishigami function, battery RUL) demonstrated the practical benefits: our exact sampling method consistently achieved higher accuracy compared to both uniform sampling and approximate leverage score techniques derived from the full tensor product matrix. This work thus provides an effective tool for substantially improving the accuracy of sketched solutions for these structured problems. While the current approach relies on the precise (permuted) lower subset structure, it opens clear avenues for future research. Promising directions include extending these efficient sampling techniques to handle more general or approximate column subset structures and automating the identification of necessary permutations to further broaden applicability.

\section{Acknowledgment}
\label{sec:acknowledgment}

This work was supported by the AFOSR awards FA9550-20-1-0138 and FA9550-20-1-0188 with Dr. Fariba Fahroo as the program manager.

\begin{appendices}

% \section{Section title of first appendix}\label{secA1}

% An appendix contains supplementary information that is not an essential part of the text itself but which may be helpful in providing a more comprehensive understanding of the research problem or it is information that is too cumbersome to be included in the body of the paper.




\end{appendices}


\bibliography{sn-bibliography}% common bib file


\begin{thebibliography}{51}
\ifx \bisbn   \undefined \def \bisbn  #1{ISBN #1}\fi
\ifx \binits  \undefined \def \binits#1{#1}\fi
\ifx \bauthor  \undefined \def \bauthor#1{#1}\fi
\ifx \batitle  \undefined \def \batitle#1{#1}\fi
\ifx \bjtitle  \undefined \def \bjtitle#1{#1}\fi
\ifx \bvolume  \undefined \def \bvolume#1{\textbf{#1}}\fi
\ifx \byear  \undefined \def \byear#1{#1}\fi
\ifx \bissue  \undefined \def \bissue#1{#1}\fi
\ifx \bfpage  \undefined \def \bfpage#1{#1}\fi
\ifx \blpage  \undefined \def \blpage #1{#1}\fi
\ifx \burl  \undefined \def \burl#1{\textsf{#1}}\fi
\ifx \doiurl  \undefined \def \doiurl#1{\url{https://doi.org/#1}}\fi
\ifx \betal  \undefined \def \betal{\textit{et al.}}\fi
\ifx \binstitute  \undefined \def \binstitute#1{#1}\fi
\ifx \binstitutionaled  \undefined \def \binstitutionaled#1{#1}\fi
\ifx \bctitle  \undefined \def \bctitle#1{#1}\fi
\ifx \beditor  \undefined \def \beditor#1{#1}\fi
\ifx \bpublisher  \undefined \def \bpublisher#1{#1}\fi
\ifx \bbtitle  \undefined \def \bbtitle#1{#1}\fi
\ifx \bedition  \undefined \def \bedition#1{#1}\fi
\ifx \bseriesno  \undefined \def \bseriesno#1{#1}\fi
\ifx \blocation  \undefined \def \blocation#1{#1}\fi
\ifx \bsertitle  \undefined \def \bsertitle#1{#1}\fi
\ifx \bsnm \undefined \def \bsnm#1{#1}\fi
\ifx \bsuffix \undefined \def \bsuffix#1{#1}\fi
\ifx \bparticle \undefined \def \bparticle#1{#1}\fi
\ifx \barticle \undefined \def \barticle#1{#1}\fi
\bibcommenthead
\ifx \bconfdate \undefined \def \bconfdate #1{#1}\fi
\ifx \botherref \undefined \def \botherref #1{#1}\fi
\ifx \url \undefined \def \url#1{\textsf{#1}}\fi
\ifx \bchapter \undefined \def \bchapter#1{#1}\fi
\ifx \bbook \undefined \def \bbook#1{#1}\fi
\ifx \bcomment \undefined \def \bcomment#1{#1}\fi
\ifx \oauthor \undefined \def \oauthor#1{#1}\fi
\ifx \citeauthoryear \undefined \def \citeauthoryear#1{#1}\fi
\ifx \endbibitem  \undefined \def \endbibitem {}\fi
\ifx \bconflocation  \undefined \def \bconflocation#1{#1}\fi
\ifx \arxivurl  \undefined \def \arxivurl#1{\textsf{#1}}\fi
\csname PreBibitemsHook\endcsname

%%% 1
\bibitem[\protect\citeauthoryear{Drineas
  et~al.}{2006}]{drineas2006SamplingAlgorithms}
\begin{bchapter}
\bauthor{\bsnm{Drineas}, \binits{P.}},
\bauthor{\bsnm{Mahoney}, \binits{M.W.}},
\bauthor{\bsnm{Muthukrishnan}, \binits{S.}}:
\bctitle{Sampling algorithms for $\ell_2$ regression and applications}.
In: \bbtitle{Proceedings of the Seventeenth Annual {{ACM}}-{{SIAM}} Symposium
  on Discrete Algorithm},
pp. \bfpage{1127}--\blpage{1136}
(\byear{2006})
\end{bchapter}
\endbibitem

%%% 2
\bibitem[\protect\citeauthoryear{Drineas
  et~al.}{2008}]{drineas2008RelativeerrorCUR}
\begin{barticle}
\bauthor{\bsnm{Drineas}, \binits{P.}},
\bauthor{\bsnm{Mahoney}, \binits{M.W.}},
\bauthor{\bsnm{Muthukrishnan}, \binits{S.}}:
\batitle{Relative-error {{CUR}} matrix decompositions}.
\bjtitle{SIAM Journal on Matrix Analysis and Applications}
\bvolume{30}(\bissue{2}),
\bfpage{844}--\blpage{881}
(\byear{2008})
\end{barticle}
\endbibitem

%%% 3
\bibitem[\protect\citeauthoryear{Drineas et~al.}{2011}]{drineas2011FasterLeast}
\begin{barticle}
\bauthor{\bsnm{Drineas}, \binits{P.}},
\bauthor{\bsnm{Mahoney}, \binits{M.W.}},
\bauthor{\bsnm{Muthukrishnan}, \binits{S.}},
\bauthor{\bsnm{Sarl{\'o}s}, \binits{T.}}:
\batitle{Faster least squares approximation}.
\bjtitle{Numerische Mathematik}
\bvolume{117}(\bissue{2}),
\bfpage{219}--\blpage{249}
(\byear{2011})
\end{barticle}
\endbibitem

%%% 4
\bibitem[\protect\citeauthoryear{Larsen and
  Kolda}{2020}]{larsen2020PracticalLeverageBased}
\begin{botherref}
\oauthor{\bsnm{Larsen}, \binits{B.W.}},
\oauthor{\bsnm{Kolda}, \binits{T.G.}}:
Practical {{Leverage}}-{{Based Sampling}} for {{Low}}-{{Rank Tensor
  Decomposition}}.
arXiv preprint arXiv:2006.16438
(2020)
{\href{https://arxiv.org/abs/2006.16438}{{arXiv:2006.16438}}}
\end{botherref}
\endbibitem

%%% 5
\bibitem[\protect\citeauthoryear{Drineas
  et~al.}{2012}]{drineas2012FastApproximation}
\begin{barticle}
\bauthor{\bsnm{Drineas}, \binits{P.}},
\bauthor{\bsnm{{Magdon-Ismail}}, \binits{M.}},
\bauthor{\bsnm{Mahoney}, \binits{M.W.}},
\bauthor{\bsnm{Woodruff}, \binits{D.P.}}:
\batitle{Fast approximation of matrix coherence and statistical leverage}.
\bjtitle{The Journal of Machine Learning Research}
\bvolume{13}(\bissue{1}),
\bfpage{3475}--\blpage{3506}
(\byear{2012})
\end{barticle}
\endbibitem

%%% 6
\bibitem[\protect\citeauthoryear{Adcock et~al.}{2022}]{adcock_sparse_2022}
\begin{bbook}
\bauthor{\bsnm{Adcock}, \binits{B.}},
\bauthor{\bsnm{Brugiapaglia}, \binits{S.}},
\bauthor{\bsnm{Webster}, \binits{C.G.}}:
\bbtitle{Sparse {Polynomial} {Approximation} of {High}-{Dimensional}
  {Functions}}.
\bpublisher{SIAM}, 
(\byear{2022})
\end{bbook}
\endbibitem

%%% 7
\bibitem[\protect\citeauthoryear{Murray
  et~al.}{2023}]{murray2023RandomizedNumerical}
\begin{botherref}
\oauthor{\bsnm{Murray}, \binits{R.}},
\oauthor{\bsnm{Demmel}, \binits{J.}},
\oauthor{\bsnm{Mahoney}, \binits{M.W.}},
\oauthor{\bsnm{Erichson}, \binits{N.B.}},
\oauthor{\bsnm{Melnichenko}, \binits{M.}},
\oauthor{\bsnm{Malik}, \binits{O.A.}},
\oauthor{\bsnm{Grigori}, \binits{L.}},
\oauthor{\bsnm{Luszczek}, \binits{P.}},
\oauthor{\bsnm{Derezi{\'n}ski}, \binits{M.}},
\oauthor{\bsnm{Lopes}, \binits{M.E.}}:
Randomized numerical linear algebra: A perspective on the field with an eye to
  software.
arXiv preprint arXiv:2302.11474v2
(2023)
{\href{https://arxiv.org/abs/2302.11474v2}{{arXiv:2302.11474v2}}}
\end{botherref}
\endbibitem

%%% 8
\bibitem[\protect\citeauthoryear{Biagioni
  et~al.}{2015}]{biagioni2015RandomizedInterpolative}
\begin{barticle}
\bauthor{\bsnm{Biagioni}, \binits{D.J.}},
\bauthor{\bsnm{Beylkin}, \binits{D.}},
\bauthor{\bsnm{Beylkin}, \binits{G.}}:
\batitle{Randomized interpolative decomposition of separated representations}.
\bjtitle{Journal of Computational Physics}
\bvolume{281}(\bissue{C}),
\bfpage{116}--\blpage{134}
(\byear{2015})
\doiurl{10.1016/j.jcp.2014.10.009}
\end{barticle}
\endbibitem

%%% 9
\bibitem[\protect\citeauthoryear{Sun et~al.}{2018}]{sun2018TensorRandom}
\begin{bchapter}
\bauthor{\bsnm{Sun}, \binits{Y.}},
\bauthor{\bsnm{Guo}, \binits{Y.}},
\bauthor{\bsnm{Tropp}, \binits{J.A.}},
\bauthor{\bsnm{Udell}, \binits{M.}}:
\bctitle{Tensor random projection for low memory dimension reduction}.
In: \bbtitle{{NeurIPS Workshop} on {Relational Representation Learning}}
(\byear{2018})
\end{bchapter}
\endbibitem

%%% 10
\bibitem[\protect\citeauthoryear{Rakhshan and
  Rabusseau}{2020}]{rakhshan2020TensorizedRandomProjectionsAISTATS}
\begin{bchapter}
\bauthor{\bsnm{Rakhshan}, \binits{B.}},
\bauthor{\bsnm{Rabusseau}, \binits{G.}}:
\bctitle{Tensorized random projections}.
In: \beditor{\bsnm{Chiappa}, \binits{S.}},
\beditor{\bsnm{Calandra}, \binits{R.}} (eds.)
\bbtitle{Proceedings of the Twenty Third International Conference on Artificial
  Intelligence and Statistics}.
\bsertitle{Proceedings of Machine Learning Research},
vol. \bseriesno{108},
pp. \bfpage{3306}--\blpage{3316}.
\bpublisher{PMLR}, 
(\byear{2020}).
\burl{https://proceedings.mlr.press/v108/rakhshan20a.html}
\end{bchapter}
\endbibitem

%%% 11
\bibitem[\protect\citeauthoryear{Rakhshan and
  Rabusseau}{2021}]{rakhshan2021RademacherRandom}
\begin{bchapter}
\bauthor{\bsnm{Rakhshan}, \binits{B.T.}},
\bauthor{\bsnm{Rabusseau}, \binits{G.}}:
\bctitle{Rademacher random projections with tensor networks}.
In: \bbtitle{{NeurIPS} Workshop on Quantum Tensor Networks in Machine Learning}
(\byear{2021})
\end{bchapter}
\endbibitem

%%% 12
\bibitem[\protect\citeauthoryear{Iwen
  et~al.}{2021}]{iwen2021LowerMemoryOblivious}
\begin{barticle}
\bauthor{\bsnm{Iwen}, \binits{M.A.}},
\bauthor{\bsnm{Needell}, \binits{D.}},
\bauthor{\bsnm{Rebrova}, \binits{E.}},
\bauthor{\bsnm{Zare}, \binits{A.}}:
\batitle{Lower memory oblivious (tensor) subspace embeddings with fewer random
  bits: Modewise methods for least squares}.
\bjtitle{SIAM Journal on Matrix Analysis and Applications}
\bvolume{42}(\bissue{1}),
\bfpage{376}--\blpage{416}
(\byear{2021})
\doiurl{10.1137/19M1308116}
{\href{https://arxiv.org/abs/https://doi.org/10.1137/19M1308116}{{https://doi.org/10.1137/19M1308116}}}
\end{barticle}
\endbibitem

%%% 13
\bibitem[\protect\citeauthoryear{Martinsson and
  Tropp}{2020}]{martinsson2020RandomizedNumerical}
\begin{botherref}
\oauthor{\bsnm{Martinsson}, \binits{P.-G.}},
\oauthor{\bsnm{Tropp}, \binits{J.}}:
Randomized numerical linear algebra: {{Foundations}} \& algorithms.
arXiv preprint arXiv:2002.01387
(2020)
{\href{https://arxiv.org/abs/2002.01387}{{arXiv:2002.01387}}}
\end{botherref}
\endbibitem

%%% 14
\bibitem[\protect\citeauthoryear{Ailon and
  Chazelle}{2009}]{ailon2009FastJohnson}
\begin{barticle}
\bauthor{\bsnm{Ailon}, \binits{N.}},
\bauthor{\bsnm{Chazelle}, \binits{B.}}:
\batitle{The fast {{Johnson}}--{{Lindenstrauss}} transform and approximate
  nearest neighbors}.
\bjtitle{SIAM Journal on Computing}
\bvolume{39}(\bissue{1}),
\bfpage{302}--\blpage{322}
(\byear{2009})
\end{barticle}
\endbibitem

%%% 15
\bibitem[\protect\citeauthoryear{Battaglino
  et~al.}{2018}]{battaglino2018PracticalRandomized}
\begin{barticle}
\bauthor{\bsnm{Battaglino}, \binits{C.}},
\bauthor{\bsnm{Ballard}, \binits{G.}},
\bauthor{\bsnm{Kolda}, \binits{T.G.}}:
\batitle{A practical randomized {{CP}} tensor decomposition}.
\bjtitle{SIAM Journal on Matrix Analysis and Applications}
\bvolume{39}(\bissue{2}),
\bfpage{876}--\blpage{901}
(\byear{2018})
\end{barticle}
\endbibitem

%%% 16
\bibitem[\protect\citeauthoryear{Jin
  et~al.}{2020}]{jin2020FasterJohnsonLindenstrauss}
\begin{barticle}
\bauthor{\bsnm{Jin}, \binits{R.}},
\bauthor{\bsnm{Kolda}, \binits{T.G.}},
\bauthor{\bsnm{Ward}, \binits{R.}}:
\batitle{Faster {{Johnson}}-{{Lindenstrauss}} transforms via kronecker
  products}.
\bjtitle{Information and Inference: A Journal of the IMA}
(\byear{2020})
\doiurl{10.1093/imaiai/iaaa028}
{\href{https://arxiv.org/abs/https://academic.oup.com/imaiai/advance-article-pdf/doi/10.1093/imaiai/iaaa028/34904656/iaaa028.pdf}{{https://academic.oup.com/imaiai/advance-article-pdf/doi/10.1093/imaiai/iaaa028/34904656/iaaa028.pdf}}}
\end{barticle}
\endbibitem

%%% 17
\bibitem[\protect\citeauthoryear{Malik and
  Becker}{2020}]{malik2020GuaranteesKronecker}
\begin{barticle}
\bauthor{\bsnm{Malik}, \binits{O.A.}},
\bauthor{\bsnm{Becker}, \binits{S.}}:
\batitle{Guarantees for the {{Kronecker}} fast
  {{Johnson}}\textendash{{Lindenstrauss}} transform using a coherence and
  sampling argument}.
\bjtitle{Linear Algebra and its Applications}
\bvolume{602},
\bfpage{120}--\blpage{137}
(\byear{2020})
\doiurl{10.1016/j.laa.2020.05.004}
\end{barticle}
\endbibitem

%%% 18
\bibitem[\protect\citeauthoryear{Bamberger
  et~al.}{2021}]{bamberger2021JohnsonLindenstraussEmbeddings}
\begin{botherref}
\oauthor{\bsnm{Bamberger}, \binits{S.}},
\oauthor{\bsnm{Krahmer}, \binits{F.}},
\oauthor{\bsnm{Ward}, \binits{R.}}:
Johnson-{L}indenstrauss embeddings with {K}ronecker structure.
arXiv preprint arXiv:2106.13349
(2021)
{\href{https://arxiv.org/abs/2106.13349}{{arXiv:2106.13349}}}
\end{botherref}
\endbibitem

%%% 19
\bibitem[\protect\citeauthoryear{Clarkson and
  Woodruff}{2017}]{clarkson2017LowRankApproximation}
\begin{barticle}
\bauthor{\bsnm{Clarkson}, \binits{K.L.}},
\bauthor{\bsnm{Woodruff}, \binits{D.P.}}:
\batitle{Low-{{Rank Approximation}} and {{Regression}} in {{Input Sparsity
  Time}}}.
\bjtitle{Journal of the ACM}
\bvolume{63}(\bissue{6}),
\bfpage{54}--\blpage{15445}
(\byear{2017})
\doiurl{10.1145/3019134}
\end{barticle}
\endbibitem

%%% 20
\bibitem[\protect\citeauthoryear{Pagh}{2013}]{pagh2013CompressedMatrix}
\begin{barticle}
\bauthor{\bsnm{Pagh}, \binits{R.}}:
\batitle{Compressed {{Matrix Multiplication}}}.
\bjtitle{ACM Transactions on Computation Theory}
\bvolume{5}(\bissue{3}),
\bfpage{9}--\blpage{1917}
(\byear{2013})
\doiurl{10.1145/2493252.2493254}
\end{barticle}
\endbibitem

%%% 21
\bibitem[\protect\citeauthoryear{Pham and Pagh}{2013}]{pham2013FastScalable}
\begin{bchapter}
\bauthor{\bsnm{Pham}, \binits{N.}},
\bauthor{\bsnm{Pagh}, \binits{R.}}:
\bctitle{Fast and {{Scalable Polynomial Kernels}} via {{Explicit Feature
  Maps}}}.
In: \bbtitle{Proceedings of the 19th {ACM SIGKDD International Conference} on
  {Knowledge Discovery} and {Data Mining}}.
\bsertitle{{{KDD}} '13},
pp. \bfpage{239}--\blpage{247}.
\bpublisher{{ACM}},
\blocation{{New York, NY, USA}}
(\byear{2013}).
\doiurl{10.1145/2487575.2487591}
\end{bchapter}
\endbibitem

%%% 22
\bibitem[\protect\citeauthoryear{Avron
  et~al.}{2014}]{avron2014SubspaceEmbeddings}
\begin{bchapter}
\bauthor{\bsnm{Avron}, \binits{H.}},
\bauthor{\bsnm{Nguyen}, \binits{H.L.}},
\bauthor{\bsnm{Woodruff}, \binits{D.P.}}:
\bctitle{Subspace {Embeddings} for the {Polynomial Kernel}}.
In: \bbtitle{Proceedings of the 27th {{International Conference}} on {{Neural
  Information Processing Systems}} - {{Volume}} 2},
pp. \bfpage{2258}--\blpage{2266}.
\bpublisher{{MIT Press}},
\blocation{{Cambridge, MA, USA}}
(\byear{2014})
\end{bchapter}
\endbibitem

%%% 23
\bibitem[\protect\citeauthoryear{Diao
  et~al.}{2018}]{diao2018SketchingKronecker}
\begin{bchapter}
\bauthor{\bsnm{Diao}, \binits{H.}},
\bauthor{\bsnm{Song}, \binits{Z.}},
\bauthor{\bsnm{Sun}, \binits{W.}},
\bauthor{\bsnm{Woodruff}, \binits{D.}}:
\bctitle{Sketching for {Kronecker Product Regression} and {P}-splines}.
In: \bbtitle{Proceedings of the 21st {International Conference} on {Artificial
  Intelligence} and {Statistics}},
pp. \bfpage{1299}--\blpage{1308}
(\byear{2018})
\end{bchapter}
\endbibitem

%%% 24
\bibitem[\protect\citeauthoryear{Ahle
  et~al.}{2020}]{ahle2020ObliviousSketchingSIAM}
\begin{bchapter}
\bauthor{\bsnm{Ahle}, \binits{T.D.}},
\bauthor{\bsnm{Kapralov}, \binits{M.}},
\bauthor{\bsnm{Knudsen}, \binits{J.B.}},
\bauthor{\bsnm{Pagh}, \binits{R.}},
\bauthor{\bsnm{Velingker}, \binits{A.}},
\bauthor{\bsnm{Woodruff}, \binits{D.P.}},
\bauthor{\bsnm{Zandieh}, \binits{A.}}:
\bctitle{Oblivious sketching of high-degree polynomial kernels}.
In: \bbtitle{Proceedings of the {Fourteenth Annual ACM-SIAM Symposium} on
  {Discrete Algorithms}},
pp. \bfpage{141}--\blpage{160}.
\bpublisher{{SIAM}}, 
(\byear{2020})
\end{bchapter}
\endbibitem

%%% 25
\bibitem[\protect\citeauthoryear{Nelson and Nguyên}{2013}]{nelson2013OSNAP}
\begin{bchapter}
\bauthor{\bsnm{Nelson}, \binits{J.}},
\bauthor{\bsnm{Nguyên}, \binits{H.L.}}:
\bctitle{Osnap: Faster numerical linear algebra algorithms via sparser subspace
  embeddings}.
In: \bbtitle{2013 IEEE 54th Annual Symposium on Foundations of Computer
  Science},
pp. \bfpage{117}--\blpage{126}
(\byear{2013}).
\doiurl{10.1109/FOCS.2013.21}
\end{bchapter}
\endbibitem

%%% 26
\bibitem[\protect\citeauthoryear{Song
  et~al.}{2021}]{song2021FastSketchingPolynomial}
\begin{bchapter}
\bauthor{\bsnm{Song}, \binits{Z.}},
\bauthor{\bsnm{Woodruff}, \binits{D.}},
\bauthor{\bsnm{Yu}, \binits{Z.}},
\bauthor{\bsnm{Zhang}, \binits{L.}}:
\bctitle{Fast sketching of polynomial kernels of polynomial degree}.
In: \bbtitle{Proceedings of the 38th International Conference on Machine
  Learning}.
\bsertitle{Proceedings of Machine Learning Research},
vol. \bseriesno{139},
pp. \bfpage{9812}--\blpage{9823}.
\bpublisher{PMLR}, 
(\byear{2021})
\end{bchapter}
\endbibitem

%%% 27
\bibitem[\protect\citeauthoryear{Ma and
  Solomonik}{2022}]{ma2022CostefficientGaussian}
\begin{botherref}
\oauthor{\bsnm{Ma}, \binits{L.}},
\oauthor{\bsnm{Solomonik}, \binits{E.}}:
Cost-efficient {G}aussian tensor network embeddings for tensor-structured
  inputs.
arXiv preprint arXiv:2205.13163
(2022)
{\href{https://arxiv.org/abs/2205.13163}{{arXiv:2205.13163}}}
\end{botherref}
\endbibitem

%%% 28
\bibitem[\protect\citeauthoryear{Diao et~al.}{2019}]{diao2019OptimalSketching}
\begin{botherref}
\oauthor{\bsnm{Diao}, \binits{H.}},
\oauthor{\bsnm{Jayaram}, \binits{R.}},
\oauthor{\bsnm{Song}, \binits{Z.}},
\oauthor{\bsnm{Sun}, \binits{W.}},
\oauthor{\bsnm{Woodruff}, \binits{D.P.}}:
Optimal {{Sketching}} for {{Kronecker Product Regression}} and {{Low Rank
  Approximation}}.
arXiv preprint arXiv:1909.13384
(2019)
{\href{https://arxiv.org/abs/1909.13384}{{arXiv:1909.13384}}}
\end{botherref}
\endbibitem

%%% 29
\bibitem[\protect\citeauthoryear{Fahrbach
  et~al.}{2022}]{fahrbach2022SubquadraticKroneckera}
\begin{botherref}
\oauthor{\bsnm{Fahrbach}, \binits{M.}},
\oauthor{\bsnm{Fu}, \binits{T.}},
\oauthor{\bsnm{Ghadiri}, \binits{M.}}:
Subquadratic {K}ronecker regression with applications to tensor decomposition.
arXiv preprint arXiv:2209.04876
(2022)
{\href{https://arxiv.org/abs/2209.04876}{{arXiv:2209.04876}}}
\end{botherref}
\endbibitem

%%% 30
\bibitem[\protect\citeauthoryear{Cheng et~al.}{2016}]{cheng2016SPALSFast}
\begin{bchapter}
\bauthor{\bsnm{Cheng}, \binits{D.}},
\bauthor{\bsnm{Peng}, \binits{R.}},
\bauthor{\bsnm{Liu}, \binits{Y.}},
\bauthor{\bsnm{Perros}, \binits{I.}}:
\bctitle{{SPALS}: {Fast} alternating least squares via implicit leverage scores
  sampling}.
In: \bbtitle{Advances {In Neural Information Processing Systems}},
pp. \bfpage{721}--\blpage{729}
(\byear{2016})
\end{bchapter}
\endbibitem

%%% 31
\bibitem[\protect\citeauthoryear{Woodruff and
  Zandieh}{2020}]{woodruff2020NearInputSparsity}
\begin{bchapter}
\bauthor{\bsnm{Woodruff}, \binits{D.}},
\bauthor{\bsnm{Zandieh}, \binits{A.}}:
\bctitle{Near input sparsity time kernel embeddings via adaptive sampling}.
In: \bbtitle{Proceedings of the 37th International Conference on Machine
  Learning}.
\bsertitle{Proceedings of Machine Learning Research},
vol. \bseriesno{119},
pp. \bfpage{10324}--\blpage{10333}.
\bpublisher{PMLR},
(\byear{2020})
\end{bchapter}
\endbibitem

%%% 32
\bibitem[\protect\citeauthoryear{Malik}{2022}]{malik2022MoreEfficientSampling}
\begin{bchapter}
\bauthor{\bsnm{Malik}, \binits{O.A.}}:
\bctitle{More efficient sampling for tensor decomposition with worst-case
  guarantees}.
In: \bbtitle{Proceedings of the 39th International Conference on Machine
  Learning}.
\bsertitle{Proceedings of Machine Learning Research},
vol. \bseriesno{162},
pp. \bfpage{14887}--\blpage{14917}.
\bpublisher{PMLR},
(\byear{2022}).
\burl{https://proceedings.mlr.press/v162/malik22a.html}
\end{bchapter}
\endbibitem

%%% 33
\bibitem[\protect\citeauthoryear{Woodruff and
  Zandieh}{2022}]{woodruff2022LeverageScoreSampling}
\begin{bchapter}
\bauthor{\bsnm{Woodruff}, \binits{D.}},
\bauthor{\bsnm{Zandieh}, \binits{A.}}:
\bctitle{Leverage score sampling for tensor product matrices in input sparsity
  time}.
In: \bbtitle{Proceedings of the 39th International Conference on Machine
  Learning}.
\bsertitle{Proceedings of Machine Learning Research},
vol. \bseriesno{162},
pp. \bfpage{23933}--\blpage{23964}.
\bpublisher{PMLR},
(\byear{2022})
\end{bchapter}
\endbibitem

%%% 34
\bibitem[\protect\citeauthoryear{Chen
  et~al.}{2020}]{PDE_inverseProblem_sketch2020}
\begin{barticle}
\bauthor{\bsnm{Chen}, \binits{K.}},
\bauthor{\bsnm{Li}, \binits{Q.}},
\bauthor{\bsnm{Newton}, \binits{K.}},
\bauthor{\bsnm{Wright}, \binits{S.J.}}:
\batitle{Structured random sketching for {PDE} inverse problems}.
\bjtitle{SIAM Journal on Matrix Analysis and Applications}
\bvolume{41}(\bissue{4}),
\bfpage{1742}--\blpage{1770}
(\byear{2020})
\doiurl{10.1137/20M1310497}
{\href{https://arxiv.org/abs/https://doi.org/10.1137/20M1310497}{{https://doi.org/10.1137/20M1310497}}}
\end{barticle}
\endbibitem

%%% 35
\bibitem[\protect\citeauthoryear{Bharadwaj et~al.}{2023}]{bharadwaj2023fast}
\begin{botherref}
\oauthor{\bsnm{Bharadwaj}, \binits{V.}},
\oauthor{\bsnm{Malik}, \binits{O.A.}},
\oauthor{\bsnm{Murray}, \binits{R.}},
\oauthor{\bsnm{Grigori}, \binits{L.}},
\oauthor{\bsnm{Buluc}, \binits{A.}},
\oauthor{\bsnm{Demmel}, \binits{J.}}:
Fast exact leverage score sampling from {K}hatri-{R}ao products with
  applications to tensor decomposition.
arXiv preprint arXiv:2301.12584v1
(2023)
{\href{https://arxiv.org/abs/2301.12584v1}{{arXiv:2301.12584v1}}}
\end{botherref}
\endbibitem

%%% 36
\bibitem[\protect\citeauthoryear{Malik and
  Becker}{2021}]{malik2021SamplingBasedMethod}
\begin{bchapter}
\bauthor{\bsnm{Malik}, \binits{O.A.}},
\bauthor{\bsnm{Becker}, \binits{S.}}:
\bctitle{A sampling-based method for tensor ring decomposition}.
In: \bbtitle{Proceedings of the 38th International Conference on Machine
  Learning}.
\bsertitle{Proceedings of Machine Learning Research},
vol. \bseriesno{139},
pp. \bfpage{7400}--\blpage{7411}.
\bpublisher{PMLR}, 
(\byear{2021}).
\burl{https://proceedings.mlr.press/v139/malik21b.html}
\end{bchapter}
\endbibitem

%%% 37
\bibitem[\protect\citeauthoryear{Meyer et~al.}{2023}]{meyer2023near}
\begin{bchapter}
\bauthor{\bsnm{Meyer}, \binits{R.A.}},
\bauthor{\bsnm{Musco}, \binits{C.}},
\bauthor{\bsnm{Musco}, \binits{C.}},
\bauthor{\bsnm{Woodruff}, \binits{D.P.}},
\bauthor{\bsnm{Zhou}, \binits{S.}}:
\bctitle{Near-linear sample complexity for lp polynomial regression}.
In: \bbtitle{Proceedings of the 2023 Annual ACM-SIAM Symposium on Discrete
  Algorithms (SODA)},
pp. \bfpage{3959}--\blpage{4025}
(\byear{2023}).
\bcomment{SIAM}
\end{bchapter}
\endbibitem

%%% 38
\bibitem[\protect\citeauthoryear{Fausett and Fulton}{1994}]{fausett_large_1994}
\begin{barticle}
\bauthor{\bsnm{Fausett}, \binits{D.W.}},
\bauthor{\bsnm{Fulton}, \binits{C.T.}}:
\batitle{Large {Least} {Squares} {Problems} {Involving} {Kronecker}
  {Products}}.
\bjtitle{SIAM Journal on Matrix Analysis and Applications}
\bvolume{15}(\bissue{1}),
\bfpage{219}--\blpage{227}
(\byear{1994})
\doiurl{10.1137/S0895479891222106}
\end{barticle}
\endbibitem

%%% 39
\bibitem[\protect\citeauthoryear{Fausett et~al.}{1997}]{fausett_improved_1997}
\begin{barticle}
\bauthor{\bsnm{Fausett}, \binits{D.W.}},
\bauthor{\bsnm{Fulton}, \binits{C.T.}},
\bauthor{\bsnm{Hashish}, \binits{H.}}:
\batitle{Improved parallel {QR} method for large least squares problems
  involving {Kronecker} products}.
\bjtitle{Journal of Computational and Applied Mathematics}
\bvolume{78}(\bissue{1}),
\bfpage{63}--\blpage{78}
(\byear{1997})
\doiurl{10.1016/S0377-0427(96)00109-4}
\end{barticle}
\endbibitem

%%% 40
\bibitem[\protect\citeauthoryear{Seshadri}{2017}]{seshadri_kronecker_2017}
\begin{botherref}
\oauthor{\bsnm{Seshadri}, \binits{P.}}:
Kronecker {Product} {Least} {Squares}.
arXiv:1705.08731 [math]
(2017).
\doiurl{10.48550/arXiv.1705.08731} .
\url{http://arxiv.org/abs/1705.08731}
\end{botherref}
\endbibitem

%%% 41
\bibitem[\protect\citeauthoryear{Marco et~al.}{2019}]{marco_least_2019}
\begin{botherref}
\oauthor{\bsnm{Marco}, \binits{A.}},
\oauthor{\bsnm{Martínez}, \binits{J.-J.}},
\oauthor{\bsnm{Viaña}, \binits{R.}}:
Least squares problems involving generalized {Kronecker} products and
  application to bivariate polynomial regression {\textbar} {SpringerLink}.
Numerical Algorithms,
21--39
(2019)
\end{botherref}
\endbibitem

%%% 42
\bibitem[\protect\citeauthoryear{Fausett and
  Hashish}{2020}]{fausett_overview_2020}
\begin{bchapter}
\bauthor{\bsnm{Fausett}, \binits{D.W.}},
\bauthor{\bsnm{Hashish}, \binits{H.}}:
\bctitle{Overview of {QR} {Methods} for {Large} {Least} {Squares} {Problems}
  {Involving} {Kronecker} {Products}}.
In: \bbtitle{Overview of {QR} {Methods} for {Large} {Least} {Squares}
  {Problems} {Involving} {Kronecker} {Products}},
pp. \bfpage{71}--\blpage{80}.
\bpublisher{De Gruyter},
(\byear{2020}).
\doiurl{10.1515/9783112314098-009}
\end{bchapter}
\endbibitem

%%% 43
\bibitem[\protect\citeauthoryear{Ghanem and
  Spanos}{2003}]{ghanem2003stochastic}
\begin{bbook}
\bauthor{\bsnm{Ghanem}, \binits{R.G.}},
\bauthor{\bsnm{Spanos}, \binits{P.D.}}:
\bbtitle{Stochastic Finite Elements: a Spectral Approach}.
\bpublisher{Courier Corporation},
(\byear{2003})
\end{bbook}
\endbibitem

%%% 44
\bibitem[\protect\citeauthoryear{Mai and Sudret}{2015}]{MaiChuVPCEF}
\begin{bchapter}
\bauthor{\bsnm{Mai}, \binits{C.V.}},
\bauthor{\bsnm{Sudret}, \binits{B.}}:
\bctitle{Polynomial chaos expansions for damped oscillators}.
In: \bbtitle{12th International Conference on Applications of Statistics and
  Probability in Civil Engineering},
\bconflocation{Vancouver, Canada}
(\byear{2015}).
\burl{https://hal.archives-ouvertes.fr/hal-01169245}
\end{bchapter}
\endbibitem

%%% 45
\bibitem[\protect\citeauthoryear{Ishigami and Homma}{1990}]{IshigamiT1990Aiqt}
\begin{bchapter}
\bauthor{\bsnm{Ishigami}, \binits{T.}},
\bauthor{\bsnm{Homma}, \binits{T.}}:
\bctitle{An importance quantification technique in uncertainty analysis for
  computer models}.
In: \bbtitle{[1990] Proceedings. First International Symposium on Uncertainty
  Modeling and Analysis},
pp. \bfpage{398}--\blpage{403}.
\bpublisher{IEEE Comput. Soc. Press}, 
(\byear{1990})
\end{bchapter}
\endbibitem

%%% 46
\bibitem[\protect\citeauthoryear{Diaz et~al.}{2018}]{Diaz2018}
\begin{barticle}
\bauthor{\bsnm{Diaz}, \binits{P.}},
\bauthor{\bsnm{Doostan}, \binits{A.}},
\bauthor{\bsnm{Hampton}, \binits{J.}}:
\batitle{Sparse polynomial chaos expansions via compressed sensing and
  d-optimal design}.
\bjtitle{Computer Methods in Applied Mechanics and Engineering}
\bvolume{336},
\bfpage{640}--\blpage{666}
(\byear{2018})
\doiurl{10.1016/j.cma.2018.03.020}
\end{barticle}
\endbibitem

%%% 47
\bibitem[\protect\citeauthoryear{Sankararaman and
  Goebel}{2013}]{Sankararaman2013UncertaintyQI}
\begin{bchapter}
\bauthor{\bsnm{Sankararaman}, \binits{S.}},
\bauthor{\bsnm{Goebel}, \binits{K.}}:
\bctitle{Uncertainty quantification in remaining useful life of aerospace
  components using state space models and inverse form}.
(\byear{2013})
\end{bchapter}
\endbibitem

%%% 48
\bibitem[\protect\citeauthoryear{Sankararaman
  et~al.}{2014}]{Sankararaman2014UQfirstOrder}
\begin{barticle}
\bauthor{\bsnm{Sankararaman}, \binits{S.}},
\bauthor{\bsnm{Daigle}, \binits{M.}},
\bauthor{\bsnm{Goebel}, \binits{K.}}:
\batitle{Uncertainty quantification in remaining useful life prediction using
  first-order reliability methods}.
\bjtitle{Reliability, IEEE Transactions on}
\bvolume{63},
\bfpage{603}--\blpage{619}
(\byear{2014})
\doiurl{10.1109/TR.2014.2313801}
\end{barticle}
\endbibitem

%%% 49
\bibitem[\protect\citeauthoryear{Sankararaman}{2015}]{Sankararaman2015significance}
\begin{barticle}
\bauthor{\bsnm{Sankararaman}, \binits{S.}}:
\batitle{Significance, interpretation, and quantification of uncertainty in
  prognostics and remaining useful life prediction}.
\bjtitle{Mechanical Systems and Signal Processing}
\bvolume{52-53},
\bfpage{228}--\blpage{247}
(\byear{2015})
\doiurl{10.1016/j.ymssp.2014.05.029}
\end{barticle}
\endbibitem

%%% 50
\bibitem[\protect\citeauthoryear{Hadigol and Doostan}{2018}]{hadigol2018least}
\begin{barticle}
\bauthor{\bsnm{Hadigol}, \binits{M.}},
\bauthor{\bsnm{Doostan}, \binits{A.}}:
\batitle{Least squares polynomial chaos expansion: A review of sampling
  strategies}.
\bjtitle{Computer Methods in Applied Mechanics and Engineering}
\bvolume{332},
\bfpage{382}--\blpage{407}
(\byear{2018})
\end{barticle}
\endbibitem

%%% 51
\bibitem[\protect\citeauthoryear{Hampton and
  Doostan}{2015}]{hampton_coherence_2015}
\begin{barticle}
\bauthor{\bsnm{Hampton}, \binits{J.}},
\bauthor{\bsnm{Doostan}, \binits{A.}}:
\batitle{Coherence motivated sampling and convergence analysis of least squares
  polynomial {Chaos} regression}.
\bjtitle{Computer Methods in Applied Mechanics and Engineering}
\bvolume{290},
\bfpage{73}--\blpage{97}
(\byear{2015})
\doiurl{10.1016/j.cma.2015.02.006}
\end{barticle}
\endbibitem

\end{thebibliography}

\end{document}